\documentclass[11pt]{amsart}

\usepackage{latexsym}
\usepackage{amscd}
\usepackage[all]{xy}
\usepackage[mathscr]{euscript}
\usepackage{dsfont}

\normalsize

\RequirePackage{xspace}
\newcommand{\thought}[1]{}
\renewcommand{\thought}[1]{ \textbf{[#1]}}

\usepackage{enumerate}        % Fancy enumerations;  used for Roman numbering.
\newenvironment{roenumerate}{\begin{enumerate}[\upshape (i)]}{\end{enumerate}}
\newcommand\nc {\newcommand}
\newcommand\rnc{\renewcommand}

\newcount\blopone
\newcount\xone
\newcount\xtwo
\newcount\ytwo

\newtheorem{theorem}{Theorem}[section]
\newtheorem{prop}[theorem]{Proposition}
\newtheorem{observation}[theorem]{Observation}
\newtheorem{com}[theorem]{Comment}
\newtheorem{redu}[theorem]{Reduction}
\newtheorem{refinement}[theorem]{Refinement}
\newtheorem{summary}[theorem]{Summary}
\newtheorem{importnota}[theorem]{Important Notation}
\newtheorem{prblm}[theorem]{Problem}
\newtheorem{notation}[theorem]{Notation}
\newtheorem{defin}[theorem]{Definition}
\newtheorem{caution}[theorem]{Caution}
\newtheorem{remark}[theorem]{Remark}
\newtheorem{reminder}[theorem]{Reminder}
\newtheorem{illustration}[theorem]{Illustration}
\newtheorem{lemma}[theorem]{Lemma}
\newtheorem{convention}[theorem]{Convention}
\newtheorem{construction}[theorem]{Construction}
\newtheorem{corollary}[theorem]{Corollary}
\newtheorem{example}[theorem]{Example}
\newtheorem{conclusion}[theorem]{Conclusion}
\newtheorem{triviality}[theorem]{Triviality}
\newtheorem{proto}[theorem]{Prototype Quasifibration}
\newtheorem{cauex}[theorem]{Cautionary Example}
\newtheorem{hypo}[theorem]{Hypothesis}
\newtheorem{subth}{ }[theorem]
\newtheorem{case}{Case}[theorem]
\newtheorem{ssubth}{ }[subth]
\newtheorem{facts}[theorem]{Facts}

\nc\tri[1]{\begin{triviality}
\label{#1}}
\nc\fac[1]{\begin{facts}
\label{#1}
\begin{em}}
\nc\cas[1]{\begin{case}
\label{#1}
\begin{em}}
\nc\cvn[1]{\begin{convention}
\label{#1}
\begin{em}}
\nc\rfn[1]{\begin{refinement}
\label{#1}}
\nc\prt[1]{\begin{proto}
\label{#1}}
\nc\lem[1]{\begin{lemma}
\label{#1}}
\nc\pro[1]{\begin{prop}
\label{#1}}
\nc\thm[1]{\begin{theorem}
\label{#1}}
\nc\obs[1]{\begin{observation}
\label{#1}}
\nc\cor[1]{\begin{corollary}
\label{#1}}
\nc\dfn[1]{\begin{defin}
\label{#1}}
\nc\sthm[1]{\begin{subth}
\label{#1}}
\nc\exm[1]{\begin{example}
\label{#1}
\begin{em}}
\nc\plm[1]{\begin{prblm}
\label{#1}
\begin{em}}
\nc\rmk[1]{\begin{remark}
\label{#1}
\begin{em}}
\nc\rmd[1]{\begin{reminder}
\label{#1}
\begin{em}}
\nc\ntn[1]{\begin{notation}
\label{#1}
\begin{em}}
\nc\smr[1]{\begin{summary}
\label{#1}
\begin{em}}
\nc\cau[1]{\begin{caution}
\label{#1}
\begin{em}}
\nc\hyp[1]{\begin{hypo}
\label{#1}}
\nc\imn[1]{\begin{importnota}
\label{#1}
\begin{em}}
\nc\rdn[1]{\begin{redu}
\label{#1}
\begin{em}}
\nc\cax[1]{\begin{cauex}
\label{#1}
\begin{em}}
\nc\cmt[1]{\begin{com}
\label{#1}
\begin{em}}
\nc\con[1]{\begin{construction}
\label{#1}
\begin{em}}
\nc\ill[1]{\begin{illustration}
\label{#1}
\begin{em}}
\nc\ssthm[1]{\begin{ssubth}
\label{#1}
\begin{em}}
\nc\cnc[1]{\begin{conclusion}
\label{#1}
\begin{em}}

\nc\elem{\end{lemma}}
\nc\erdn{\end{em}\end{redu}}
\nc\erfn{\end{refinement}}
\nc\eprt{\end{proto}}
\nc\ethm{\end{theorem}}
\nc\eobs{\end{observation}}
\nc\ecor{\end{corollary}}
\nc\edfn{\end{defin}}
\nc\esthm{\end{subth}}
\nc\epro{\end{prop}}
\nc\etri{\end{triviality}}
\nc\eexm{\end{em}
\end{example}}
\nc\ecvn{\end{em}
\end{convention}}
\nc\ecmt{\end{em}
\end{com}}
\nc\efac{\end{em}
\end{facts}}
\nc\ermk{\end{em}
\end{remark}}
\nc\ermd{\end{em}
\end{reminder}}
\nc\eill{\end{em}
\end{illustration}}
\nc\eplm{\end{em}
\end{prblm}}
\nc\ecas{\end{em}
\end{case}}
\nc\ecau{\end{em}
\end{caution}}
\nc\ecax{\end{em}
\end{cauex}}
\nc\eimn{\end{em}
\end{importnota}}
\nc\entn{\end{em}
\end{notation}}
\nc\econ{\end{em}
\end{construction}}
\nc\esmr{\end{em}
\end{summary}}
\nc\ehyp{
\end{hypo}}
\nc\ecnc{\end{em}
\end{conclusion}}
\nc\essthm{\end{em}
\end{ssubth}}

\nc\sst{\scriptstyle}
\newcommand{\comment}[1]{}
\newcommand{\ri}{\longrightarrow}

\newcommand{\zz}{{\mathbb Z}}
\nc\gm{{{\mathbb G}_m}}

\newcommand{\D}{{\mathbf D}}
\newcommand{\Dqc}{{\mathbf D_{\text{\bf qc}}}}

\newcommand{\pp}{{\mathbb P}}

\newcommand{\oo}{\otimes}

\nc\op{^{\hbox{\rm\tiny op}}}
\nc\mth{^{\hbox{\rm\tiny th}}}

\nc\dcoh{\D_{\mathbf{coh}}^b}
\nc\script{\mathscr}
\nc\z{\zeta}
\nc\bc{{\mathbb{BC}}}
\nc\ct{{\script T}}
\nc\cf{{\script F}}
\nc\cl{{\script L}}
\nc\cv{{\script V}}
\nc\cq{{\script Q}}
\nc\cu{{\script U}}
\nc\ce{{\script E}}
\nc\cg{{\script G}}
\nc\ch{{\script H}}
\nc\cs{{\script S}}
\nc\car{{\script R}}
\nc\cd{{\script D}}
\nc\cc{{\script C}}
\nc\ck{{\script K}}
\nc\ca{{\script A}}
\nc\ci{{\script I}}
\nc\cj{{\script J}}
\nc\co{{\script O}}
\nc\cm{{\script M}}
\nc\cz{{\script Z}}
\nc\bd{\begin{description}}
\nc\ed{\end{description}}
\nc\ctob{{\script C}at\big(\ci^{op},\ca\big)}
\nc\clim{{\ds\mathop{\rm lim}_{\ds\longleftarrow}}\,}
\nc\climi{\clim^{\!i}\,}
\nc\climn{\clim^{\!n}\,}
\nc\colim{{\ds\mathop{\rm colim}_{\ds\la}}}
\nc\oa{\overline{\ca}}
\nc\s{\sigma}
\nc\ta{\tau}
\nc\os{\overline\sigma}
\nc\ot{\overline\tau}
\nc\T{\Sigma}
\nc\Tm{\Sigma^{-1}}
\nc\de[1]{{\mathop{\rm deg(#1)}}}
\nc\Ad[1]{\mathop{\rm Ad}(#1)}
\nc\ad[1]{\mathop{\rm ad}(#1)}
\nc\kth{{\it K}--theory}
\nc\loc[1]{{\text{\rm Loc(#1)}}}
\nc\coloc[1]{{\text{\rm Coloc}(#1)}}

\nc\one{\mathds{1}}

\def\der #1 {D\left(#1\right)}
\nc\prf{\begin{proof}}
\nc\eprf{\end{proof}}
\nc\ds{\displaystyle}
\nc\Tor{\text{\rm Tor}}

\nc\cb{{\script B}}
\nc\ab{{\script A}b}

\nc\be{\begin{roenumerate}}
\nc\ee{\end{roenumerate}}

\nc\cat[1]{{\script C}at\Big({\big\{#1\big\}}\op\,\,,\,\,\ab\Big)}
\nc\csab{{\script C}at\big(\cs^{op},\ab\big)}
\nc\ctab{{\script C}at\Big({\{\ct^\alpha\}}^{op},\ab\Big)}
\nc\csex{{\script E}x\big(\cs^{op},\ab\big)}
\nc\ctex{{\script E}x\Big({\{\ct^\alpha\}}^{op},\ab\Big)}
\nc\sub{\qquad\subset\qquad}
\nc\ctr[1]{{\left.\ct\left(-,#1\right)\right|}_{\cs}}
\nc\ctrf[2]{{\left.\ct\left(#1,#2\right)\right|}_{\cs}}
\nc\Ctr[1]{{\left.\ct\left(-,#1\right)\right|}_{\ct^\alpha}}
\nc\Ctrf[2]{{\left.\ct\left(#1,#2\right)\right|}_{\ct^\alpha}}

\nc\la{\longrightarrow}
\nc\nin{\noindent}
\nc\cad[1]{\text{card}(#1)}
\nc\eq{\quad=\quad}
\nc\BA{\begin{array}{c}}
\nc\EA{\end{array}}
\nc\barr{
\[
\begin{array}{cccccccccccccccc}
}
\nc\earr{
\end{array}
\]
}
\nc\as[1]{{\langle S\rangle}^{#1}}
\nc\sh{\text{\it shift}}

\nc\yy[1]{{\left.\ct\left(-,#1\right)\right|}_{\ct^c}}
\nc\vrep[2]{{\left.\ct\left(#1,#2\right)\right|}_{\ct^\alpha}}
\nc\da{\downarrow}
\nc\Hom{{\mathop{\rm Hom}}}
\nc\RHom{{\mathop{\rm RHom}}}
\nc\HHom{{\script H}{\mathop{\rm om}}}
\nc\RHHom{{\script{RH}}{\mathop{\rm om}}}
\nc\End{{\mathop{\rm End}}}
\nc\Ext{{\mathop{\rm Ext}}}
\nc\mr\Modtc
\nc\PExt{{\mathop{\rm PExt}}}
\nc\stm{\text{\rm stmod}(kG)}
\nc\stM{\text{\rm StMod}(kG)}
\nc\e{\varepsilon}
\nc\p{\mathfrak{p}}

\nc\rs{\s^{-1}A}
\nc\br{{\{\s^{-1}A\}}}
\nc\ra\ri
\nc\y[1]{\mathbf{y}#1}
\nc\x[1]{\mathbf{z}#1}
\nc\mmod[1]{#1\text{--\rm mod}}
\nc\Mod[1]{#1\text{--\rm Mod}}
\nc\Md {\ensuremath{\mathop{\textup{Mod}}}}
\rnc\mod[1]{\ensuremath{\mathop{\textup{mod-}#1}}\xspace}
\nc\Modtc{\Mod{\ct^c}}
\nc\pgldim[1]{\mathop{\rm pgldim}\,#1}
\nc\tf{{\rm [TR5]}}
\nc\tfs{{\rm [TR5$^*$]}}
\nc\Fun{\text{\rm Funct}(F\op,\ab)}
\nc\sym{\text{\rm Sym}}
\nc\sgn{\text{\rm sgn}}
\nc\Pro{\text{\rm Prod}^{}_\alpha(F\op,\ab)}
\nc\Yt[1]{{\left.\Hom_\ct^{}\left(-,#1\right)\right|}_F^{}}
\nc\dl{\delta}
\nc\Proj[1]{#1\text{--\rm Proj}}
\nc\proj[1]{#1\text{--\rm proj}}
\nc\Flat[1]{\text{\rm Flat}\,#1}
\nc\Inj[1]{\text{\rm Inj\,}#1}
\nc\qc[1]{\text{\rm qc\,}#1}
\nc\ov{\overline}
\nc\wt{\widetilde}
\nc\ph{\varphi}
\nc\tstr{{\it t}--structure}
\nc\spec[1]{{\text{\rm Spec}(#1)}}
\nc\Spec[1]{{\text{\rm Spec}\big(#1\big)}}

\newcommand{\Dqcpl}{{\mathbf D^+_{\mathbf{qc}}}}
\newcommand{\Dqcmi}{{\mathbf D^-_{\mathbf{qc}}}}

\nc\EProd{\text{\rm EProd}}
\nc\ECoprod{\text{\rm ECoprod}}
\nc\Prod{\text{\rm Prod}}
\nc\Coprod{\text{\rm Coprod}}
\nc\ldimp{\text{\rm LDim}^{\prod}}
\nc\ldimc{\text{\rm LDim}^{\coprod}}

\nc\hoco{
\begin{picture}(40,10)
\put(20,0){\makebox(0,0)[b]{\text{\rm Hocolim}}}
\put(5,-2){\vector(1,0){30}}
\end{picture}\,}

\nc\holim{
\begin{picture}(40,10)
\put(20,0){\makebox(0,0)[b]{\text{\rm Holim}}}
\put(35,-2){\vector(-1,0){30}}
\end{picture}}

\nc\Cop{\text{\rm Coprod}}
\nc\seq{{\mathbb{S}_{\mathbf{e}}}}
\nc\se{{S^{\text{\tt e}}}}
\nc\te{{T^{\text{\tt e}}}}
\nc\LL{{\text{\bf L}}}
\nc\R{\text{\bf R}}
\nc\id{\text{\rm id}}
\nc\ev{\text{\rm ev}}
\nc\supp{\text{\rm supp}}
\nc\Loc{\text{\rm Loc}}
\nc\Thick{\text{\rm Thick}}

\begin{document}

\author{Amnon Neeman}\thanks{The research was partly supported 
by the Australian Research Council}
\address{Centre for Mathematics and its Applications \\
        Mathematical Sciences Institute\\
        John Dedman Building\\
        The Australian National University\\
        Canberra, ACT 0200\\
        AUSTRALIA}
\email{Amnon.Neeman@anu.edu.au}

\title[Grothendieck duality and Hochschild homology]{The relation 
            between Grothendieck 
            duality and Hochschild homology}

\begin{abstract}
Grothendieck duality goes back to 1958, to the talk at the 
ICM in Edinburgh~\cite{Grothendieck58C} announcing the result.
Hochschild homology is even older,
its roots can be traced back to the 1945 article~\cite{Hochschild45}.
The fact that the two might be related is relatively recent.
The first hint of a relationship came in 1987 in
Lipman~\cite{Lipman87}, 
and another was found in 1997 in Van den 
Bergh~\cite{vandenBergh97}. Each of these discoveries was 
interesting and had an impact, Lipman's mostly
by giving another approach to the computations
and Van den Bergh's especially
on the development 
of non-commutative versions of the subject.  However in this survey
we will almost entirely focus on a third, much more recent connection,
discovered in 2008
by Avramov and Iyengar~\cite{Avramov-Iyengar08}
and later developed and extended in several papers,
see for example \cite{Avramov-Iyengar-Lipman-Nayak10,Iyengar-Lipman-Neeman13}.

There are two classical paths to the foundations of
Grothendieck duality, one following
Grothendieck and Hartshorne~\cite{Hartshorne66} and (much later) 
Conrad~\cite{Conrad00}, and the other following Deligne~\cite{Deligne66},
Verdier~\cite{Verdier68} and (much later) Lipman~\cite{Lipman09}. The
accepted view is that each of these has its 
drawbacks: the first approach (of Grothendieck, Hartshorne and Conrad) is
complicated and messy to set up, while the second (of Deligne, Verdier
and Lipman) might be cleaner to present but leads to a theory where it's not
obvious how to compute anything.

The point of this article is that the recently-discovered connection
with Hochschild homology and cohomology (the one due to Avramov
and Iyengar) changes this. It renders
clearly superior the highbrow appoach to the
subject, the one due to Deligne, Verdier and Lipman. Not only
is it (relatively) easy to set up the machinery, the computations also
become transparent. And in the process we learn that Grothendieck
duality is not really about residues of
meromorphic differential forms, it is about
the local cohomology of the Hochschild homology. By a fortuitous accident,
if $f:X\la Y$ is a smooth map then the top Hochschild homology happens
to be isomorphic to the relative canonical bundle, and its top local cohomology 
is represented by meromorphic differential forms. This is the reason
that, as long as we stick to smooth maps, what comes up
is residues of meromorphic
forms. For non-smooth, flat maps it's Hochschild homology and maps from
it that we need to study. 
\end{abstract}



\maketitle

\tableofcontents

\setcounter{section}{-1}

\section{Introduction}
\label{S-1}

Let $f:X\la Y$ be a morphism of noetherian schemes. At the 
level of derived categories 
there exist
 natural functors $\LL f^*:\Dqc(Y)\la\Dqc(X)$, its right adjoint 
$\R f_*:\Dqc(X)\la
\Dqc(Y)$, 
as well as a right adjoint for $\R f_*$, nowadays (following Lipman)
denoted 
$f^\times:\Dqc(Y)\la\Dqc(X)$. 
For general $f$ the functor
$f^\times$ can be dreadful---it can take
a bounded complex of coherent sheaves, that is an
object in $\dcoh(Y)\subset\Dqc(Y)$,
to a truly enormous object in $\Dqc(X)$.
This functor $f^\times$ only behaves well under
strong restrictions, the usual being that $f$ be
proper.

To remedy this one introduces a better-behaved functor $f^!$. If $f$ is
proper then $f^!=f^\times$, but for general $f$ one traditionally
does some finicky
manipulations to arrive at $f^!$. And, until very recently, the recipe
worked only for cohomologically bounded-below complexes. That is
$f^!$ has always been viewed as a functor $f^!:\Dqcpl(Y)\la\Dqcpl(X)$.

Against this background came the striking work of 
Avramov, Iyengar, Lipman 
and Nayak, see~\cite{Avramov-Iyengar08,Avramov-Iyengar-Lipman-Nayak10},
relating Grothendieck duality with Hochschild homology
and cohomology. To give the flavor of
the results let me present just one formula, and for simplicity
let me give only the affine version.
Suppose therefore that $X=\spec S$, $Y=\spec R$, assume that $R$ and
$S$ are noetherian, and that $f:X\la Y$ is a flat, finite-type
map. In an abuse of notation we will write $f:R\la S$ for the induced
ring homomorphism, and also identify $\D(R)\cong\Dqc(Y)$ and $\D(S)=\Dqc(X)$.
Let $\se=S\oo_R^{}S$ be the enveloping algebra. Then, for any
object $N\in\Dqcpl(Y)=\D^+(R)$, we have a canonical isomorphism
\[
f^!N\cong S\oo_{\se}^{}\Hom_R^{}(S,S\oo_R^{}N)\ .
\]
In this formula the tensor products and the Hom are all derived.

The reader might find it interesting to note that, in the special case
where $f:R\la S$ is finite and \'etale, we recover the classical
formula
\[
f^!N\cong \Hom_R^{}(S,N)\cong S\oo_{\se}^{}\Hom_R^{}(S,S\oo_R^{}N)\ .
\]
Of course for finite, \'etale maps the Homs and tensors are underived.
We will revisit \'etale maps (not necessarily finite) in Remark~\ref{R1.13}.

Perhaps one needs some familiarity with the classical literature to
appreciate how striking this is---assuming only that $f$ is
flat we  have produced a formula for $f^!$, which took
a mere paragraph to state, and is clearly free of
auxiliary choices and
functorial. And although the left-hand-side was defined on the
assumption that $N$ is bounded
below---after all we only knew $f^!$ on the
bounded-below derived category---the right-hand-side makes sense for any $N$.
In fact the formula tells us the surprising fact that if $N$ is an
object in $\D^+(R)$ then $S\oo_{\se}^{}\Hom_R^{}(S,S\oo_R^{}N)$ must belong
to $\D^+(S)$. We know, from the complicated
classical construction, that $f^!$ takes
$\D^+(R)$ to $\D^+(S)$, but $S$ is not of finite Tor-dimension
over $\se$ and we have
no reason to expect an expression of the form $S\oo_{\se}^{}M$ to be
bounded below. The derived tensor product tends to introduce lots of
negative cohomology.

In joint work with Iyengar and Lipman we revisited these results,
and along the way developed a useful
new natural transformation $\psi(f):f^\times\la f^!$, 
see~\cite{Iyengar-Lipman-Neeman13}. Hints of $\psi$ may be found in 
Lipman~\cite[Exercise~4.2.3(d)]{Lipman09}, but without the naturality
properties that make it so valuable.
With all these unexpected new tools
it was becoming clear that the time may have come to revisit
the foundations of Grothendieck duality.
In this article we sketch what has come out of this.

Finally we should tell the reader the structure of this survey. The
early sections, \S\ref{S1} and \S\ref{S2}, survey recent
results that can be found elsewhere in the literature.  The results
are new, meaning new in this generality---there are older avatars, 
what's unusual here
is that the theory is developed in the unbounded derived category. The
results might be innovative but we still omit the proofs.
With the
exception of Proposition~\ref{P2.202}, where the argument 
is included, the proofs are all to be found
in recent preprints available electronically.

In \S\ref{S99} and \S\ref{S95} this changes. Special cases of the
results are known, with what turn out to be artificial boundedness
restrictions.
We give a general treatment---both to show that the
results are true more generally, and to illustrate the power of the
new techniques. Because of this our treatment is complete, with
proofs.
The reader interested in the highlights is advised to read the
statement (not proof) of Lemma~\ref{L99.5}, as well as
Corollary~\ref{C95.13} and Example~\ref{E95.15}. 

The final sections, \S\ref{S3} and \S\ref{S4}, are again ``soft'',
with no proofs presented. They review the history and suggest open problems.

\section{Conventions}
\label{S-2}

In this article we consider schemes $X$ and the corresponding
derived categories $\Dqc(X)$, whose objects are complexes
of sheaves of $\co_X^{}$--modules
with quasicoherent cohomology. Since abelian categories never come
up, whenever there is a possible ambiguity our functors should be
assumed derived---thus we will write $f^*$ for $\LL f^*$, $f_*$
for $\R f_*$,  $\Hom$ for $\RHom$ and $\oo$ for the derived tensor
product $\oo^\LL_{}$.
For simplicity, in \S\ref{S1}, \S\ref{S2}, \S\ref{S99}, \S\ref{S95}
and \S\ref{S3} we
will assume that our schemes are noetherian
and morphisms of schemes are separated and of finite
type---occasionally, but not always, we will explicitly remind the reader of
these standing assumptions.
Unless we specifically say otherwise
all derived categories will be unbounded. For a 
morphism of schemes $f:X\la Y$ we let 
$f^*\dashv f_*\dashv f^\times$ be the adjoint functors
which, back in \S\ref{S-1}, we referred to as
$\LL f^*\dashv \R f_*\dashv f^\times$.
 
\section{The formal theory}
\label{S1}

In this section and the next
we sketch the current state of
the formal
theory, without worrying about who proved 
what and when. 

Let $f:X\la Y$ be a morphism of schemes. The functor $f^*:\Dqc(Y)\la\Dqc(X)$
is a strict monoidal functor, meaning it respects the tensor product. 
Therefore for any pair of objects $E\in\Dqc(Y)$ and $F\in\Dqc(X)$ we have a
natural map
\[
\CD
f^*[E\oo f_*F] @>\sim>> f^*E\oo f^*f_*F @>\id\oo\e>> f^*E\oo F 
\endCD
\] 
where the first map is the natural isomorphism, and $\e:f^*f_*\la\id$ is 
the counit of the adjunction $f^*\dashv f_*$. By adjunction we obtain a
natural map $p(E,F):E\oo f_*F\la f_*(f^*E\oo F)$. The map $p(E,F)$
is known to be an
isomorphism, usually called the
\emph{projection formula.} This leads us to

\dfn{D1.1}
Let $f:X\la Y$ be a morphism of schemes and let $E,F$ be objects in
$\Dqc(Y)$. The map $\chi(f,E,F):f^*E\oo f^\times F\la f^\times (E\oo F)$
is defined by applying the adjunction $f_*\dashv f^\times$ to the 
composite
\[
\xymatrix@C+30pt{ 
f_*(f^*E\oo f^\times F) \ar[r]^-{p(E,f^\times F)^{-1}} & E\oo f_*f^\times F
\ar[r]^-{\id\oo\e'} & E\oo F\ ,
}
\] 
where the first map is the inverse of the isomorphism in the projection
formula, while $\e':f_*f^\times\la \id$ is the counit of the
adjunction $f_*\dashv f^\times$.
\edfn

The first result in the theory is 

\thm{T1.3}
The map $\chi(f,E,F)$ is an isomorphism whenever
\be
\item
$f$ is arbitrary, but $E$ is a perfect complex.
\item
$E$ and $F$ are arbitrary, but $f$ is proper and of finite Tor-dimension.
\ee
\ethm

Next recall the base-change maps. Given a  commutative square
of morphisms of schemes
\[
\CD
W @>u>> X \\
@VfVV  @VVgV \\
Y @>v>> Z 
\endCD
\]
there is a canonical isomorphism of functors $\alpha:f^*v^*\la u^*g^*$. Consider
the composite
\[
\CD
f^*v^*g_* @>\alpha g_*>> u^*g^*g_* @>u^*\e>> u^*\ , 
\endCD
\]
where $\e:g^*g_*\la\id$ is the counit of the adjunction $g^*\dashv g_*$. 
Adjunction gives us a base-change map $\beta:v^*g_*\la f_*u^*$; this map
is not always an isomorphism, but there are important situations
in which it is. This leads us to

\dfn{D1.5}
Assume we are in a situation where the base-change map $\beta:v^*g_*\la f_*u^*$
is an isomorphism; for this article the important case 
where this happens is when the square
\[
\CD
W @>u>> X \\
@VfVV  @VVgV \\
Y @>v>> Z 
\endCD
\]
is cartesian and the map $v$ is flat. In this scenario consider the 
composite
\[
\xymatrix@C+20pt{
f_*u^*g^\times \ar[r]^-{\beta^{-1}g^\times} & v^*g_*g^\times \ar[r]^-{v^*\e'} &
v^*
}
\]
where the first map is the inverse of the isomorphism $\beta$ while
$\e':g_*g^\times\la\id$ is the counit of the adjunction
$g_*\dashv g^\times$.
The (second) base change map $\Phi:u^*g^\times\la f^\times v^*$ corresponds
to this composite under the adjunction $f_*\dashv f^\times$.
\edfn

One can wonder when the base-change map $\Phi$ is an isomorphism. 
The best result to date says

\thm{T1.7}
Let the notation be as in the case of  
Definition~\ref{D1.5} which interests us in this article---that 
is we assume the square cartesian and $v$ flat.
Let $E$ be an object
in $\Dqc(Z)$. Then the
 base-change map $\Phi(E):u^*g^\times(E)\la f^\times v^*(E)$
is an isomorphism provided $g$ is proper
and one of the conditions below holds:
\be
\item
$E$ belongs to $\Dqcpl(Z)\subset\Dqc(Z)$.
\item
$E\in\Dqc(Z)$ is arbitrary, but the map $f:W\la Y$ is of finite Tor-dimension.
\ee 
\ethm

Now one proceeds as follows:
given any morphism $f:X\la Y$ we factor it as 
$X\stackrel u\la\ov X\stackrel p\la Y$ with $u$ an open immersion
and $p$ proper, and then define $f^!:\Dqc(Y)\la\Dqc(X)$ by
the formula $f^!=u^*p^\times$. One of the consequences of 
 Theorem~\ref{T1.7} is that $f^!$
is well-defined, meaning that it is canonically independent of
the choice of factorization. And we have the following
theorem.

\thm{T1.107}
The assignment, taking a morphism of schemes $f:X\la Y$
to the functor $f^!:\Dqc(Y)\la\Dqc(X)$, satisfies a long list
of compatibility properties. We list some highlights.
\sthm{ST1.107.1}
Let 
$X\stackrel f\la Y\stackrel g\la Z$
be composable morphisms of schemes.
There is a map $\rho(f,g):{(gf)}^!\la f^!g^!$, which has the property
that the two ways of using $\rho$ to go
from $(hgf)^!$ to $f^!g^!h^!$ are equal.
\esthm
\sthm{ST1.107.3}
The two functors $f^\times, f^!:\Dqc(Y)\la\Dqc(X)$ are related
by a natural transformation $\psi(f):f^\times\la f^!$. The $\psi$ is compatible
with composition, in the obvious sense that the square below commutes
\[\xymatrix@C+15pt{
{(gf)}^\times \ar[r]^{\delta(f,g)}\ar[d]_{\psi(gf)} &f^\times g^\times \ar[d]^{\psi(f)\psi(g)}\\
{(gf)}^! \ar[r]_{\rho(f,g)} & f^!g^!
}\]
where $\rho(f,g)$ is the map of \ref{ST1.107.1}, while $\delta(f,g):(gf)^\times
\la f^\times g^\times$ is the canonical isomorphism.
\esthm
\sthm{ST1.107.4}
The map $\rho(f,g)$ is an
isomorphism if $f$ is of finite Tor-dimension or if either
$gf$ or $g$ is proper. The map $\psi(f)$ is an isomorphism
whenever $f$ is proper.
\esthm
\sthm{ST1.107.5}
Given a pair of object $E,F\in\Dqc(Y)$ then there is a way to mimick
the construction in Definition~\ref{D1.1} with $f^!$ in place of $f^\times$.
More precisely: there is a map $\s(f,E,F):f^*E\oo f^!F\la f^!(E\oo F)$ so that
the natural square commutes
\[
\CD
f^*E\oo f^\times F @>\chi(f,E,F)>> f^\times(E\oo F)\\
@V\id\oo\psi(f)VV  @VV\psi(f)V \\
 f^*E\oo f^!F @>\sigma(f,E,F)>> f^!(E\oo F)
\endCD
\]
Furthermore we have the analog of Theorem~\ref{T1.3}, that is
$\s(f,E,F)$ is an isomorphism if one of the conditions below
holds
\be
\item
$f$ is arbitrary, but $E$ is a 
perfect complex.
\item
$E$ and $F$ are arbitrary, but 
$f$ is of finite Tor-dimension. 
\ee
\esthm
\sthm{ST1.107.7}
The base-change map $\Phi$ of Definition~\ref{D1.5} also has an $(-)^!$
analog. Precisely: given a cartesian square as in Definition~\ref{D1.5},
there is a base-change map $\theta:u^*g^!\la f^!v^*$. 
\esthm
\sthm{ST1.107.9}
There is an analog of Theorem~\ref{T1.7} for $(-)^!$ in place of $(-)^\times$.
Precisely: the map $\theta(E):u^*g^!(E)\la f^!v^*(E)$ is an isomorphism
as long as one of the following holds
\be
\item
$E$ belongs to $\Dqcpl(Z)\subset\Dqc(Z)$.
\item
$E\in\Dqc(Z)$ is arbitrary, but the map $f:W\la Y$ is of finite Tor-dimension.
\ee 
\esthm
\ethm

The full list of compatibility properties is quite long, and in any case it
is clearer and more compact to present it in a 2-category formulation.
For this paper we content ourselves with what's in Theorem~\ref{T1.107}. 

\rmk{R1.109}
In the introduction we mentioned that people have traditionally preferred
$f^!$ to $f^\times$ because it is ``better behaved''. Theorem~\ref{T1.107}
allows us to make this more precise. If we compare \ref{ST1.107.5} with
Theorem~\ref{T1.3} we see that
\be
\item
  If $f$ is proper then $\s(f,E,F)$ and $\chi(f,E,F)$ agree up to canonical
  isomorphism.
  To see this observe that, when $f$ is proper, then the vertical
  maps in the commutative square of \ref{ST1.107.5} are isomorphisms
  by \ref{ST1.107.4}.
\item
  The maps $\s(f,E,F)$ and $\chi(f,E,F)$ are defined for every
  triple $f,E,F$, but
  $\s(f,E,F)$ is an isomorphism more often. If $E$ is perfect
  then both are isomorphisms. But for non-perfect $E$ the result
  \ref{ST1.107.5}(ii)
  says that $\s(f,E,F)$ is
  an isomorphism whenever $f$ is of finite
  Tor-dimension, whereas Theorem~\ref{T1.3}(ii)
  guarantees that $\chi(f,E,F)$ is an isomorphism only if $f$ is
  proper as well as of finite Tor dimension.
\setcounter{enumiii}{\value{enumi}}
\ee
The same pattern repeats itself for the base-change maps $\Phi$ and $\theta$.
They are defined for every cartesian square with flat horizontal morphisms,
and coincide if the vertical
maps are proper---if our Theorem~\ref{T1.107} were less pared down this could
be shown to follow from the general structure, the reader can see the
introduction to \cite{Neeman13} for the fullblown formalism.
But if we ask ourselves when $\Phi$ and $\theta$ induce isomorphisms,
the conditions on
$\theta$ are less restrictive than on $\Phi$. Precisely: when we compare
Theorem~\ref{T1.7} with  \ref{ST1.107.9} we discover
\be
\setcounter{enumi}{\value{enumiii}}
\item
  Assume $E$ is bounded below. Then 
  $\Phi(E):u^*g^\times E\la f^\times v^*E$ is an isomorphism if $g$ is proper,
  while $\theta(E):u^*g^! E\la f^! v^*E$ is an isomorphism unconditionally.
\item
  Let $E$ be arbitrary. Then 
  $\Phi(E)$ is an isomorphism as long as
  $g$ is proper and $f$ is of
  finite Tor-dimension, while $\theta(E)$ is an isomorphism whenever
  $f$ is of
  finite Tor-dimension (no need for any properness).
\setcounter{enumiii}{\value{enumi}}
\ee
\ermk

\rmk{R1.11}
If $f:X\la Y$ is an open immersion then the square 
\[
\xymatrix@C+2pt@R+2pt{
X \ar[r]^{\id}\ar[d]_{\id} & X\ar[d]^{f}\\
X \ar[r]^{f} & Y
}\]
is cartesian. By \ref{ST1.107.9}(ii) we have that 
$f^!=\id^*f^!\stackrel{\theta}\la \id^! f^*=f^*$ 
is an isomorphism.
Given any morphism $g:Y\la Z$, the map 
$\rho(f,g):(gf)^!\la f^!g^!$ is an isomorphism by
\ref{ST1.107.4}---after all the open immersion
$f$ is of finite Tor-dimension.
Combining these isomorphisms gives
\be
\item
  If $X\stackrel f\la Y\stackrel g\la Z$ are composable
  morphisms of schemes, with $f$ an open immersion, then
  we have a canonical isomorphism $(gf)^!\stackrel{\rho}\la
  f^!g^!\stackrel\theta\la f^*g^!$.
\ee
If $g:Y\la Z$ happens to be a proper morphism
then $\psi(g):g^\times\la g^!$ is an isomorphism by \ref{ST1.107.4},
which we may combine with the isomorphism of (i) to deduce
a canonical isomorphism
 $(gf)^!\cong f^*g^!\cong f^*g^\times$. The compatibilities
of Theorem~\ref{T1.107} force upon us the formula for $(gf)^!$. That is:
any time we can factor a map $X\la Z$ as a composite
$X\stackrel f\la Y\stackrel g\la Z$, with $f$ an open immersion and $g$
proper, then $(gf)^!$ must be given by the formula $(gf)^!\cong f^*g^\times$.
\ermk

\rmk{R1.13}
In passing we observe that the
formula of Remark~\ref{R1.11}(i) generalizes, we need only
assume $f$ \'etale.
Suppose $f:X\la Y$ is an \'etale morphism of noetherian schemes. Consider the
following diagram 
\[
\xymatrix@C+15pt@R+2pt{
X\ar[r]^-{\Delta} & X \times_Y^{}X\ar[r]^-{\pi_1^{}}\ar[d]_-{\pi_2^{}} & X\ar[d]^{f}\\
 &X \ar[r]^{f} & Y
}\]
where the square is cartesian and $\Delta$ is the diagonal map. We are
given that $f$ is \'etale, meaning that the diagonal map $\Delta$ is
an open immersion. This gives a series of isomorphisms
\[
f^! \,\,\cong\,\, \Delta^*\pi_1^*f^! \,\,\cong\,\, \Delta^*\pi_2^!f^* \,\,\cong\,\,
(\pi_2^{}\Delta)^! f^* \,\,\cong\,\, f^*.
\]
The first isomorphism is $\Delta^*\pi_1^*\cong\id^*=\id$, the second
isomorphism is the map $\pi_1^*f^!\la 
\pi_2^!f^*$  of \ref{ST1.107.7}, which is an isomorphism by \ref{ST1.107.9}(ii), the third
isomorphism is Remark~\ref{R1.11}(i) applied to the composable maps
$X\stackrel\Delta\la X\times_Y^{} X\stackrel{\pi_2^{}}\la X$ with
$\Delta$ an open immersion, and the last isomorphism is because
$\pi_2^{}\Delta$ is the identity. Remark~\ref{R1.11}(i) therefore also
generalizes, we have
\be
\item
  If $X\stackrel f\la Y\stackrel g\la Z$ are composable morphisms with
  $f$ \'etale, then $(gf)^!\stackrel\rho\la f^!g^!\cong f^*g^!$
  combines to an isomorphism. The map $\rho(f,g)$ is an isomorphism
  because $f$ is of finite Tor-dimension, while
  $f^!\cong f^*$ is the isomorphism above. 
\ee
\ermk

\section{Still formal, but less familiar---the way the abstract theory is related to explicit computations}
\label{S2}

Our first aim is to obtain a fomula for $f^!$ free
of auxiliary choices---that is, one that does not involve
factoring $f$ as $X\stackrel u\la\ov X\stackrel p\la Y$
with $u$ an open immersion and $p$ proper. We begin with a little Lemma.

\lem{L2.2}
Let $U\stackrel\alpha\la V\stackrel\beta\la W$ be finite-type morphisms of
noetherian schemes
so that $\alpha$ is a closed immersion and $\beta\alpha$ is
proper. Then the maps $\alpha^\times\psi(\beta):\alpha^\times\beta^\times
\la\alpha^\times\beta^!$ and  
$\alpha^*\psi(\beta):\alpha^*\beta^\times
\la\alpha^*\beta^!$ 
are both isomorphisms.
\elem

\prf
The proof is an easy consequence of Theorem~\ref{T1.107}, coupled
with standard facts from support theory. 
First~\ref{ST1.107.3} gives us the commutative square
\[\xymatrix@C+15pt{
{(\beta\alpha)}^\times \ar[r]^{\delta(\alpha,\beta)}\ar[d]_{\psi(\beta\alpha)} &\alpha^\times \beta^\times \ar[d]^{\psi(\alpha)\psi(\beta)}\\ 
{(\beta\alpha)}^! \ar[r]_{\rho(\alpha,\beta)} & \alpha^!\beta^!
}\]
Recalling
that $\alpha$ and $\beta\alpha$ are proper maps,
we conclude from \ref{ST1.107.4}
that $\rho(\alpha,\beta):(\beta\alpha)^!\la\alpha^!\beta^!$ is an isomorphism,
as are $\psi(\alpha):\alpha^\times\la\alpha^!$ and
$\psi(\beta\alpha):(\beta\alpha)^\times\la(\beta\alpha)^!$.
Now $\delta(\alpha,\beta)$
is trivially an isomorphism, hence in the square above
the indicated maps are all isomorphisms
\[\xymatrix@C+15pt{
{(\beta\alpha)}^\times \ar[r]^{\sim}\ar[d]_{\wr} &\alpha^\times \beta^\times \ar[d]^{\psi(\alpha)\psi(\beta)}\\ 
{(\beta\alpha)}^! \ar[r]_{\sim} & \alpha^!\beta^!
}\]
The commutativity implies that the vertical map 
$\psi(\alpha)\psi(\beta):\alpha^\times\beta^\times
\la
\alpha^!\beta^!$ 
must be an isomorphism. This isomorphism can be written as the
composite 
\[
\CD
\alpha^\times\beta^\times
@>{\alpha^\times\psi(\beta)}>>
\alpha^\times\beta^!
@>{\psi(\alpha)\beta^!}>>
\alpha^!\beta^!\ .
\endCD
\]
but, as $\psi(\alpha)$ is an isomorphism, 
we conclude that $\alpha^\times\psi(\beta)$
is also an isomorphism. Support theory, more precisely 
\cite[Proposition~A.3(ii)]{Iyengar-Lipman-Neeman13}, tells us that
$\alpha^*\psi(\beta)$ is also an isomorphism. 
\eprf

We will apply this little lemma in the following situation.

\con{C2.2002}
Let $f:X\la Y$ be a finite-type, flat morphism of noetherian schemes.
We may form the diagram
\[\xymatrix@C+10pt@R+10pt{
X \ar[dr]^-{\Delta} & & \\
 & X\times_Y^{}X \ar[r]^-{\pi_1^{}}\ar[d]_-{\pi_2^{}} & X\ar[d]^f \\
 & X \ar[r]^-{f} & Y\ar@{}[ul]|{(\diamondsuit)}
}\]
where the square is cartesian, $\pi_1^{}$ and $\pi_2^{}$ are the first
and second projections, and $\Delta:X\la X\times_Y^{}X$ is the
diagonal inclusion. We assert
\econ

\pro{P2.202}
With the notation as in Construction~\ref{C2.2002},
there is a canonical isomorphism $\Delta^*\pi_2^\times f^*\la f^!$.
\epro

\prf
We apply Lemma~\ref{L2.2} to the
composable maps $X\stackrel\Delta\la  X\times_Y^{}X\stackrel {\pi_2^{}}\la X$;
the map $\Delta$ is a closed immersion and the composite $\id=\pi_2^{}\Delta$
is proper, and Lemma~\ref{L2.2} tells us that
$\Delta^*\psi(\pi_2^{}):\Delta^*\pi_2^\times\la \Delta^*\pi_2^!$ is an isomorphism.

On the other hand all the maps in the cartesian square $(\diamondsuit)$
are flat, and 
\ref{ST1.107.9}(ii) gives that
$\theta:\pi_1^*f^!\la\pi_2^!f^*$ is an isomorphism.
Consider therefore the composite
\[
\xymatrix@C+25pt@R+5pt{
  \Delta^*\pi_2^\times f^* \ar[r]^-{\Delta^*\psi(\pi_2^{})f^*} &
  \Delta^*\pi_2^! f^* \ar[r]^-{\Delta^*\theta^{-1}_{}} &
  \Delta^*\pi_1^* f^! \ar[r]^-{\sim} & \id^* f^!\ar@{=}[r] & f^! 
}
\]
The first and second maps are isomorphisms by the discussion above,
and the
third map comes from applying the functor $(-)^*$ to the equality
$\id=\pi_1^{}\Delta$. The composite is therefore
an isomorphism $\Delta^*\pi_2^\times f^*\la f^!$.
\eprf

\con{C2.3}
Assume $f:X\la Y$ is flat and let the notation be as above.
If $\e:\Delta_*\Delta^\times\la\id$
is the counit of adjunction $\Delta_*\dashv\Delta^\times$, let 
$\ph:\Delta_*\la\pi_2^\times$ be the composite 
\[
\xymatrix@C+20pt{
\Delta_*  \ar@{=}[r] & \Delta_*\Delta^\times\pi_2^\times 
\ar[r]^-{\e\pi_2^\times} & \pi_2^\times
}
\]
where the equality is the observation that 
$\id=\id^\times=(\pi_2^{}\Delta)^\times
\cong\Delta^\times\pi_2^\times$. 
Define $c_f^{}:\Delta^*\Delta_*f^*\la\Delta^*\pi_2^\times f^*$
to be $\Delta^*\ph f^*$; combining with the isomorphism
$\Delta^*\pi_2^\times f^*\la f^!$ of Proposition~\ref{P2.202} we
arrive at a map which, in an abuse
of notation, we will also denote $c_f^{}:\Delta^*\Delta_*f^*\la f^!$.
If we apply this map to the object $\co_Y^{}\in\Dqc(Y)$ and note
that $f^*\co_Y^{}=\co_X^{}$, we obtain a map
 $c_f^{}(\co_Y^{}):\Delta^*\Delta_*\co_X^{}\la f^!\co_Y^{}$.
For any integer $d\in\zz$ let $\gamma_f^{}(d)$ be the composite
\[
\xymatrix@C+20pt{
(\Delta^*\Delta_*\co_X^{})^{\leq -d}\ar[r] &
\Delta^*\Delta_*\co_X^{}\ar[r]^-{c_f^{}(\co_Y^{})} & f^!\co_Y^{}
}
\]
where the first map is given by the standard \tstr\ truncation.
\econ

Note that we have defined the maps in great generality, globally and
without auxiliary choices of coordinates. What is known so far is

\thm{T2.5}
Suppose the map $f:X\la Y$ is smooth and of relative dimension $d$.
Then the map $\gamma_f^{}(d)$ is an isomorphism.
\ethm

\rmk{R2.7}
The object $f^!\co_Y^{}$ might appear mysterious but $\Delta^*\Delta_*\co_X^{}$
isn't, it is just the obvious object in the derived category whose cohomology
is the Hochschild homology of $\co_X^{}$. If $f:X\la Y$ is smooth and of
relative dimension $d$ then $(\Delta^*\Delta_*\co_X^{})^{\leq -d}$ is 
nothing other than $HH^d(\co_X^{})$, which is the relative canonical bundle
in degree $-d$. In symbols 
$\Omega^d_{X/Y}[-d]=(\Delta^*\Delta_*\co_X^{})^{\leq -d}$.

Note that the maps $c_f^{}$ and $\gamma_f^{}(d)$ are defined for any flat 
$f$ and any integer $d$, and might contain interesting information for
$f$ which aren't smooth. I don't believe anyone has computed examples yet.
\ermk

\rmk{R2.9}
The reader can find (different) proofs of Theorem~\ref{T2.5} in either
Alonso, Jerem{\'{\i}}as and
Lipman~\cite[Proposition~2.4.2]{Alonso-Jeremias-Lipman14} or else in
\cite[\S1]{Neeman13A}.
The point we want to make here is that the proof can't be hard:
the map is defined globally, but proving it an isomorphism
is local in $Y$ in the flat topology---hence we may assume $Y$ an
affine scheme---and 
 local in $X$ in the \'etale topology\footnote{To see that it suffices
  to prove the map an isomorphism flat-locally in $Y$ and \'etale-locally
  in $X$ one needs the
  full strength of Theorem~\protect{\ref{T1.107}}, the extract we presented
  in this survey doesn't suffice. The reader can find the complete statement
  in~\protect{\cite{Neeman13}}}.
And \'etale-locally
any smooth map of degree $d$ is of the
form $\spec S\la \spec R$, where
$f:R\la S$ identifies $S$ as the polynomial ring 
$S=R[x_1^{},x_2^{},\ldots,x_d^{}]$.
Let $\se=S\oo_R S$; since $S$ is flat over $R$ it makes
no difference whether we view this particular
tensor product as ordinary or derived.
The expression $\Delta^*\Delta_*\co_X^{}$ is nothing other than
the derived tensor product $S\oo_{\se}^{}S$, while 
$f^!\co_Y^{}\cong\Delta^*\pi_2^\times f^*\co_Y^{}$ comes down to
$S\oo_{\se}^{}\Hom_R^{}(S,S)$, where the tensor and the Hom are
both derived. And the map $c_f^{}:S\oo_{\se}^{}S\la
S\oo_{\se}^{}\Hom_R^{}(S,S)$ 
is the tensor product over $\se$ of the identity
map $\id:S\la S$ and the obvious inclusion
$S\la \Hom_R^{}(S,S)$. 
OK: a little computation is necessary to finish off the proof,
the details may be found
in \cite[\S1]{Neeman13A}.
\ermk

\con{C2.11}
Now let $W\subset X$ be the union of closed subsets $W_i\subset X$, 
such that
the restriction of $f$ to each $W_i$ is proper. 
For every point $x\in X$ write $k(x)$ for its residue field;
then the full subcategory $\D_{\mathrm{qc},W}^{}(X)\subset\Dqc(X)$
has objects given by the formula 
\[
\D_{\mathrm{qc},W}^{}(X)=\{E\in\Dqc(X)\mid E\oo k(x)=0\text{ for all }x\notin W\}.
\]
The inclusion
 $I:\D_{\mathrm{qc},W}^{}(X)\la\Dqc(X)$ is well-known to admit a Bousfield
localization, meaning it has a right adjoint 
$R:\Dqc(X)\la\D_{\mathrm{qc},W}^{}(X)$. The composite functor
$\Dqc(X)\stackrel R\la\D_{\mathrm{qc},W}^{}(X)\stackrel I\la\Dqc(X)$ 
is nowadays denoted
$\Gamma_W^{}:\Dqc(X)\la \Dqc(X)$,\footnote{ 
Classically it was denoted $\R\Gamma_W^{}$---it is the right derived functor
of some functor on abelian categories.
But we have been suppressing all the notation
that usually reminds us of the functors of abelian categories that
we are deriving, and in this case it brings our notation into
concert with that of Benson, Iyengar and
Krause~\cite{Benson-Iyengar-Krause08,Benson-Iyengar-Krause11,Benson-Iyengar-Krause11A,Benson-Iyengar-Krause12,Benson-Iyengar-Krause13}. Their
choice of the letter $\Gamma$ was quite unrelated to Grothendieck's,
see the comment at the top of
\cite[page~582]{Benson-Iyengar-Krause08}.
By a fortuitous accident the notations coincide
(once the $\R$ is eliminated in $\R\Gamma$).}
and this Bousfield localization is
even smashing, meaning there is a natural isomorphism
$E\oo\Gamma_W^{}F\la\Gamma_W^{}(E\oo F)$. We assumed that, for each
$W_i$, the composite $W_i\stackrel{\alpha_i}\la X\stackrel f\la Y$
is proper, and
a slight refinement of Lemma~\ref{L2.2} tells us that $\alpha_i^*\psi(f)$ is
an isomorphism for each $\alpha_i$. Support theory, more precisely
\cite[Proposition~A.3(ii)]{Iyengar-Lipman-Neeman13}, allows us to
deduce that
the map $\Gamma_W^{}\psi(f)$ is also an isomorphism.

Let $\e:\Gamma_W^{}=IR\la\id$ be the counit of the adjunction $I\dashv R$
and let $\e':f_*f^\times\la\id$ be the counit of the adjunction 
$f_*\dashv f^\times$. We define the map $\int_W^{}:f_*\Gamma_W^{}f^!\la\id$ to
be the composite
\[
\xymatrix@C+40pt{
f_*\Gamma_W^{}f^!\ar[r]_-{f_*\big[\Gamma_W^{}\psi(f)\big]^{-1}}
\ar@/^2pc/[rrr]^-{\ds\int_W^{}}
 & :f_*\Gamma_W^{}f^\times \ar[r]_-{f_*\e f^\times} & f_*f^\times
\ar[r]_-{\e'} &\id
}\]
If we apply this natural transformation to the
object $\co_Y^{}\in\Dqc(Y)$, and combine with the map
$\gamma_f^{}(d)$, we obtain a composite
\[
\xymatrix@C+40pt{
f_*\Gamma_W^{}\big[\Delta^*\Delta_*\co_X^{}\big]^{\leq -d}
\ar[r]^-{f_*\Gamma_W^{}\gamma_f^{}(d)} &
f_*\Gamma_W^{}f^!\co_Y^{}
\ar[r]^-{\ds\int_W^{}}
&\co_Y^{}
}\]
Now assume $f:X\la Y$ is smooth of relative dimension $d$. Theorem~\ref{T2.5}
tells us that $\gamma_f^{}(d)$ is an isomorphism,
and in Remark~\ref{R2.7} we observed that 
$\big[\Delta^*\Delta_*\co_X^{}\big]^{\leq -d}$ is nothing more
than the relative canonical bundle in degree $-d$, that
is $\Omega^d_{X/Y}[-d]=\big[\Delta^*\Delta_*\co_X^{}\big]^{\leq -d}$.
We are assuming $W=\cup W_i$ is the union 
of closed subschemes $W_i\subset X$ so that, for each $i$, the composite
$W_i\stackrel{\alpha_i}\la X\stackrel f\la Y$ is proper. Let us make the
stronger assumption that $f\alpha_i$ is a finite morphism for each $i$. Then
the functor $f_*\Gamma_W^{}$ takes an object $E\in\Dqc(X)$ to its
local cohomology at $W$ and, in the particular case of $E=\Omega^d_{X/Y}[-d]$,
this comes down (locally and modulo boundaries)
 to the relative meromorphic $d$--forms
$\omega/s_1^{}s_2^{}\cdots s_d^{}$ on $X$ with $(s_1^{},s_2^{},\ldots, s_d^{})$
a relative system of parameters. And now we are ready for the next computation.
\econ

\thm{T2.13}
Assume $f:X\la Y$ is smooth of relative dimension $d$.
Then the composite
\[
\xymatrix@C+0pt{
f_*\Gamma_W^{}\Omega^d_{X/Y}[-d]\ar@{=}[r] &
f_*\Gamma_W^{}\big[\Delta^*\Delta_*\co_X^{}\big]^{\leq -d}
\ar[rr]^-{f_*\Gamma_W^{}\gamma_f^{}(d)} & &
f_*\Gamma_W^{}f^!\co_Y^{}
\ar[rr]^-{\ds\int_W^{}}&
&\co_Y^{}
}\]
is just the map taking a relative meromorphic $d$--form to its residue.
\ethm

\rmk{R2.15}
Early incarnations of Theorem~\ref{T2.13}  may be found in
Verdier~\cite[pp.~398-400]{Verdier68} and H\"ubl and
Sastry~\cite[Residue Theorem, p.~752]{Hubl-Sastry93}.
The proof of Theorem~\ref{T2.13}, as stated here, may be found in \cite[\S2]{Neeman13A}. Once
again we note that while the global definition might be slightly subtle,
the local computation can easily be reduced to the simple, affine case of 
Remark~\ref{R2.9}.
\ermk

We have mentioned, several times, that the functor $f^!$ is ``better behaved''
than its cousin $f^\times$. One facet is that it is amenable
to local computations. To illustrate this we end the section with
a couple of little lemmas. The first of the lemmas just records, for a
morphism $f:X\la Y$, the open subsets $U\subset X$ which we will find
useful for these local computations.

\lem{L2.2055}
Suppose $f:X\la Y$ is a finite-type morphism of
noetherian schemes. 
Suppose $U\subset X$ is an open affine subset so that the composite
$U\stackrel u\la X\stackrel f\la Y$ can be factored through an open affine subset $V\subset Y$.
Then
\be
\item
  If $f$ is of finite Tor-dimension then $fu$ admits a factorization
  $fu=hg$, where $g$ is a closed immersion of finite Tor-dimension and
  $h$ is smooth.
\item
  For arbitrary $f$, the map $fu$ may be factored as $fu=khg$, where $g$
  is an open immersion, $h$ is a closed immersion and $k$ is smooth and
  proper.
\ee
\elem

\prf
By hypothesis the composite $U\stackrel u\la X\stackrel f\la Y$ is
equal to the composite 
$U\stackrel \alpha\la V\stackrel\beta\la Y$ with $V$ affine and $\beta$
an open immersion.
Because $\alpha:U\la V$ is a finite-type morphism of affine schemes we may
factor it
as $U\stackrel {j'}\la \pp^n_V\stackrel{\pi'}\la V$ where
$j'$ is a locally closed
immersion, and this allows us to factor $fu$ as
$U\stackrel {j}\la \pp^n_Y\stackrel{k}\la Y$
where $j$ is also a locally closed immersion.

Now we separate the treatment into two cases. If $f$ is of finite
Tor-dimension we factor $j$ as
$U\stackrel g\la W\stackrel \delta\la \pp_Y^{n}$ with $g$ a
closed immersion and $\delta$ an open immersion. This gives
us a factorization of $fu$ as $U\stackrel g\la W\stackrel{k\delta}\la Y$.
Since $k$ and $\delta$ are both smooth so is $h=k\delta$. And
the fact that $fu=hg$ is of finite Tor-dimension and $h$ is smooth means
that $g$ must be of finite Tor-dimension. We have found a factorization
satisfying (i) of the Lemma.

It remains to treat the case where $f$ is arbitrary. In
this case 
we factor
$j:U\la\pp_Y^n$ as $U\stackrel g\la W\stackrel h\la \pp_Y^{n}$, with
$g$ an open immersion and $h$ a closed immersion.
In total this factors $fu$ as $khg$, with $k:\pp_Y^n\la Y$
smooth and proper, $h:W\la\pp_Y^n$ a closed immersion and
$g:U\la W$ an open immersion.
\eprf

\lem{L1.13}
Let $f:X\la Y$ be a finite-type, separated morphism of noetherian
schemes. We record the following boundedness and coherence statements:
\be
\item
  There exists an integer $n$ so that $f^!\Dqc(Y)^{\geq 0}\subset\Dqc(X)^{\geq n}$.
\item
  If $E\in\Dqc(Y)$ is bounded below and has coherent cohomology, then the
  same is true for $f^!E$.
\setcounter{enumiii}{\value{enumi}}
\ee
For the next few assertions assume furthermore that $f$ is
of finite Tor-dimension.
\be
\setcounter{enumi}{\value{enumiii}}
\item
  For any object $E\in\Dqc(Y)$, the support of $f^!E$ is equal to
  the support of $f^*E$.
\item
  If $E\in\Dqc(Y)$ has coherent cohomology then so has
  $f^!E$ [no need to assume $E$
  bounded below.]
\item
  There exists an integer $n$ so that $f^!\Dqc(Y)^{\leq 0}\subset\Dqc(X)^{\leq n}$.
\item
  There exists an integer $n$ so that, if the Tor-amplitude of $E\in\Dqc(Y)$
  is contained in the interval $[0,\infty)$, then the
    Tor-amplitude of $f_*f^!E$ is contained in the interval
    $[n,\infty)$.
\setcounter{enumiii}{\value{enumi}}      
\ee
\elem

\prf
First we show that (vi) follows from (i) and Theorem~\ref{T1.107}.
Suppose we know (i);  we may choose an integer $n$ with
$f^!\Dqc(Y)^{\geq0}\subset\Dqc(X)^{\geq n}$. Let $E\in\Dqc(Y)$
be an object with Tor-amplitude contained in
$[0,\infty)$. Then $\Dqc(Y)^{\geq0}\oo E\subset\Dqc(Y)^{\geq0}$,
and
\[
f^*\Dqc(Y)^{\geq0}\oo f^!E\eq f^!\big[\Dqc(Y)^{\geq0}\oo E\big]\quad\subset\quad
\Dqc(X)^{\geq n}
\]
where the equality is by \ref{ST1.107.5}(ii).
Applying $f_*$ we deduce
\[
\Dqc(Y)^{\geq0}\oo f_*f^!E\eq
f_*\big[f^*\Dqc(Y)^{\geq0}\oo f^!E\big] \quad\subset\quad
f_*\Dqc(X)^{\geq n}\quad\subset\quad \Dqc(Y)^{\geq n}
\]
where the equality is by the projection formula. We conclude that
$f_*f^!E$ has Tor-amplitude in the interval $[n,\infty)$.

It remains to prove (i)--(v), which are all local in $X$.
This means the following: cover $X$ by open immersions $u_i:U_i\subset X$. 
Remark~\ref{R1.11}(i) tells us that $(fu_i)^!\cong u_i^*f^!$.
If we had a cover so that we could prove (i)--(v) for all
the $(fu_i)^!\cong u_i^*f^!$, then the statement would follow for $f$.
Cover $Y$ by finitely many open affine subsets $V_i$, then cover each
$f^{-1}V_i$ by finitely many open affine subsets $U_i$, and we have covered $X$ by open
subsets as in Lemma~\ref{L2.2055}. We are therefore
reduced to proving the Lemma under the assumption that
$f=fu$ has a factorization as in Lemma~\ref{L2.2055}~(i) or (ii).

We first prove (i), (ii), (iv) and (v), all of which
respect composition: this means
\be
\setcounter{enumi}{\value{enumiii}}
\item
  If $X\stackrel f\la Y\stackrel g\la Z$
are morphisms so that 
the map $\rho(f,g):(gf)^!\la f^!g^!$ is an isomorphism
and 
(i) and (ii)  hold for each of $f$ and $g$,
then it formally follows that (i) and (ii)  hold for $gf$.
\item
If we furthermore assume that $f$ and $g$ are of finite Tor-dimension
and (iv) and (v) hold for each of $f$ and $g$, then (iv) and (v)
hold for $gf$.
\setcounter{enumiii}{\value{enumi}}      
\ee
Next we observe that the factorizations of Lemma~\ref{L2.2055}
behave well with respect to the functor $(-)^!$, meaning
\be
\setcounter{enumi}{\value{enumiii}}
\item
  If $f=hg$ is a factorization as in Lemma~\ref{L2.2055}(i)
  then 
  $\rho(g,h):(hg)^!\la g^!h^!$ is an isomorphism. This is because
  $g$ is of finite Tor-dimension, see \ref{ST1.107.4}.
\item
If $f=khg$ is a factorization as in Lemma~\ref{L2.2055}(ii)
  then 
  the map $(khg)^!\la g^!h^!k^!$ is an isomorphism.
  The map $\rho(hg, k):(khg)^!\la(hg)^!k^!$ is
  an isomorphism because $k$ is proper,
  while the map $\rho(g,h):(hg)^!\la g^!h^!$ is an isomorphism
  because $g$ is of finite Tor-dimension. See \ref{ST1.107.4}.
\setcounter{enumiii}{\value{enumi}}      
\ee
This means that it suffices to prove (i), (ii), (iv) and (v)
in the special cases where $f$ is either smooth or a closed immersion,
and for the proof of (iv) and (v) we may assume that the closed
immersion is of finite Tor-dimension.

If $f:X\la Y$ is a smooth map
then \ref{ST1.107.5}(ii) tells us that $f^!(-)\cong f^*(-)\oo f^!\co_Y^{}$,
while from Theorem~\ref{T2.5} we learn that $f^!\co_Y^{}$ is just
a shift of the relative
canonical bundle. Hence (i), (ii), (iv) and (v) are all true
for smooth maps $f$.

Next we prove (i), (ii), (iv) and (v) for closed immersions $f$. By
\ref{ST1.107.4} we know that $\psi(f):f^\times\la f^!$
is an isomorphism. We need to prove the coherence and/or vanishing
of cohomology sheaves of $f^!E\cong f^\times E$, and
as $f$ is a closed immersion these are equivalent to the
coherence and/or vanishing of the cohomology sheaves
of $f_*f^\times E$.
Now recall the isomorphisms
$f_*f^\times(-)\cong f_*\HHom\big[\co_X^{},f^\times(-)\big]
\cong\HHom(f_*\co_X^{},-)$.
If $f$ is of finite Tor-dimension then $f_*\co_X^{}\in\Dqc(Y)$ is a perfect
complex, hence (iv) and (v) are clear. It remains to prove
(i) and (ii), which become assertions about the vanishing and coherence
of $\HHom(f_*\co_X^{},E)$ where $E$ is bounded below. 
Illusie~\cite[Proposition~3.7]{Illusie71A} tells us that this
may be computed locally in $Y$, and if $Y$ is affine
the assertions are obvious.

It remains to prove (iii), the assertion about the supports.
The map $f$ is assumed of finite Tor-dimension and \ref{ST1.107.5}(ii)
gives an isomorphism $f^!(E)\cong f^*(E)\oo f^!\co_Y^{}$.
Support theory tells us that the support of the tensor
product $f^*(E)\oo f^!\co_Y^{}$ is the intersection of the support
of $f^*E$ and the support of $f^!\co_Y^{}$, hence it suffices
to show that the support of $f^!\co_Y^{}$ is all of $X$.
But we have reduced to the case where 
$f$ has
a factorization $X\stackrel g\la W\stackrel h\la Y$ as
in Lemma~\ref{L2.2055}(i), and in (ix) we noted that
 $\rho(g,h):(hg)^!\la g^!h^!$ is an
isomorphism. Thus
\[
(hg)^!\co_Y^{}\,\cong\, g^![h^!\co_Y^{}]\,\cong\, \big[g^*h^!\co_Y^{}\big]\oo g^!\co_W^{}
\]
where the last isomorphism is by \ref{ST1.107.5}(ii). We wish to
show that the support of $(hg)^!\co_Y^{}$ is all of $X$,
and it suffices to prove that the support of
$h^!\co_Y^{}$ is all of $W$ and the support of $g^!\co_W^{}$ is
all of $X$.

Theorem~\ref{T2.5} tells us that $h^!\co_Y^{}$ is just
a shift of the relative canonical bundle---its support is all of $W$.
Since $g$ is a closed immersion the support of
$g^!\co_W^{}\cong g^\times\co_W^{}$ is equal
to the support of $g_*g^\times\co_W^{}\cong\HHom(g_*\co_X^{},\co_W^{})$. But
as $g$ has finite Tor-dimension the object $g_*\co_X^{}$ is a perfect
complex on $W$, and its support (all of $X$) is equal to the support of the
dual $\big[g_*\co_X^{}\big]^\vee=\HHom(g_*\co_X^{},\co_W^{})$. This completes
the proof of (iii).
\eprf

\section{Some basic isomorphisms}
\label{S99}

In this section we prove some formal corollaries of the theory
presented so far, establishing that certain natural maps are
isomorphisms. None of the results is hard, but they 
are a little technical---their value will only become
apparent when we see the applications later on. The one result we
will need, in \S\ref{S95}, is Lemma~\ref{L99.5}. At a first
reading we recommend that the reader study the statement of 
Lemma~\ref{L99.5} and skip the rest of this section.

\lem{L99.-1}
Let $f:X\la Y$ be a finite-type morphism of noetherian schemes, of
finite Tor-dimension. Let
$E\in\Dqc(Y)$ be a perfect complex, and let $\{A_\lambda,\,\lambda\in\Lambda\}$
be a set of objects of $\Dqc(Y)$. Then
the natural map
$\Big[\prod_{\lambda\in\Lambda}^{}f^*A_\lambda\Big]\oo f^!E
\la \prod_{\lambda\in\Lambda}^{}\big[f^*A_\lambda\oo f^!E\big]$
is an isomorphism.
\elem

\prf
We begin by proving a special case: we show that the Lemma is true
$f$ is proper and of finite Tor-dimension.
Assume therefore that $f$ is proper and of finite Tor-dimension. Then $f_*$
takes perfect complexes to perfect complexes,
and \cite[Theorem~1.7(GN2)]{Balmer-Dellambrogio-Sanders16} tells us that $f^*$
has a left adjoint and respects products.
If we contemplate the commutative diagram
\[\xymatrix@C+20pt{
 \ds f^*\left[\prod_{\lambda\in\Lambda}A_\lambda\right]\oo f^! E
  \ar[r]^-{(1)}\ar[d]_{\s}&
  \ds\left[\prod_{\lambda\in\Lambda}f^*A_\lambda\right]\oo f^! E
  \ar[r]^-{(2)}&
  \ds\prod_{\lambda\in\Lambda}\big[f^*A_\lambda\oo f^! E\big]\ar[d]^\s\\
  \ds f^!\left(\left[\prod_{\lambda\in\Lambda}A_\lambda\right]\oo E\right)
  \ar[r]^-{(3)} &
  \ds f^!\left[\prod_{\lambda\in\Lambda}\big(A_\lambda\oo E\big)\right]
  \ar[r]^-{(4)} &
  \ds\prod_{\lambda\in\Lambda}f^!\big[A_\lambda\oo E\big]
}\]
the vertical maps $\s$ are isomorphisms by \ref{ST1.107.5}(ii), the map
$(1)$ is an isomorphism because $f^*$ 
respects products, the map $(3)$ is an isomorphism because $E$
is perfect and hence $(-)\oo E\cong\HHom(E^\vee,-)$ respects products,
and the map $(4)$ is an isomorphism because, for the proper morphism
$f$, the functor $f^!\cong f^\times$ has a left adjoint and respects
products. Hence the map $(2)$ must be an isomorphism, as we asserted.

We have proved the Lemma if $f$ is proper
and of finite Tor-dimension. Now suppose $f$ may be factored
as $X\stackrel g\la W\stackrel h\la Y$
with $h$ smooth and $g$ proper and of finite Tor-dimension.
Because $h$ is of finite Tor-dimension \ref{ST1.107.5} gives an
isomorphism $h^*E\oo h^!\co_Y^{}\la h^!E$, while
Theorem~\ref{T2.5} tells us that $h^!\co_Y^{}$ is a shift
of the relative
canonical bundle. Since $E$ is assumed perfect
it follows that $h^!E$, being the
tensor product of two perfect complexes $h^*E$ and
$h^!\co_Y^{}$, must be perfect. But then we may apply
the
special case
of the Lemma to the map $g:X\la W$, the perfect complex $h^!E\in\Dqc(W)$
and the set of objects $\{h^*A_\lambda,\,\lambda\in\Lambda\}$
of $\Dqc(W)$. We deduce that the map $(2)$ in the
commutative square below is an isomorphism
\[\xymatrix@C+30pt{
  \ds\left(\prod_{\lambda\in\Lambda}\big[(hg)^*A_\lambda\big]\right)\oo (hg)^!E
  \ar[d]_{\tau(g,h)\oo\rho(g,h)}
  \ar[r]^-{(1)} &\ds\prod_{\lambda\in\Lambda}\big[(hg)^*A_\lambda\oo (hg)^!E\big]
  \ar[d]^{\tau(g,h)\oo\rho(g,h)}\\
\ds\left(\prod_{\lambda\in\Lambda}\big[g^*h^*A_\lambda\big]\right)\oo g^!h^!E
\ar[r]^-{(2)} &\ds\prod_{\lambda\in\Lambda}\big[g^*h^*A_\lambda\oo g^!h^!E\big]
}\]
The vertical maps are induced by the isomorphism $\tau(g,h):(hg)^*\la g^*h^*$
and the map $\rho(g,h):(hg)^!\la g^!h^!$ of \ref{ST1.107.1}.
Because $g$ is of finite Tor-dimension \ref{ST1.107.4}
informs us that $\rho(g,h)$ is an isomorphism, hence
both vertical maps are isomorphisms.
From the commutativity we deduce that $(1)$ is an isomorphism.
In other words: the Lemma is true for any $f$ which admits a factorization
as $X\stackrel g\la W\stackrel h\la Y$ with $h$ smooth and
$g$ proper and of finite Tor-dimension.

Now let $f:X\la Y$ be arbitrary, fix a perfect
complex $E\in\Dqc(Y)$, as well as
a set of objects $\{A_\lambda,\,\lambda\in\Lambda\}$
of $\Dqc(Y)$. If $u:U\la X$ is an open immersion
with $U$ affine, and assuming that $fu:U\la Y$ can be factored through
an open affine subset $V\subset Y$, then Lemma~\ref{L2.2055}(i)
guarantees that $fu$ may be factored as $hg$ with $h$ smooth and
$g$ proper and of finite Tor-dimension. By the above the natural
map
$\big[\prod_{\lambda\in\Lambda}(fu)^*A_\lambda\big]\oo(fu)^!E
  \la\prod_{\lambda\in\Lambda}\big[(fu)^*A_\lambda\oo(fu)^!E\big]$
is an 
isomorphism. Remark~\ref{R1.11}(i) gives an
isomorphism $(fu)^!\cong u^*f^!$, allowing us to rewrite the  
isomorphism above as
$\big[\prod_{\lambda\in\Lambda}(fu)^*A_\lambda\big]\oo u^*f^!E
\la\prod_{\lambda\in\Lambda}\big[(fu)^*A_\lambda\oo u^*f^!E\big]$.
If we apply the functor $u_*$, and use the projection formula
and the fact that $u_*$ has a left adjoint and hence
respects products, we obtain that the natural map
is an isomorphism
$\big[\prod_{\lambda\in\Lambda}u_*u^*f^*A_\lambda\big]\oo f^!E
\la\prod_{\lambda\in\Lambda}\big[u_*u^*f^*A_\lambda\oo f^!E\big]$.

Now glue these isomorphisms: any time we have open
immersions $u:U\la X$, $v:V\la X$, $j:U\cap V\la X$ and $w:U\cup V\la X$,
then for each $\lambda$ we obtain a triangle
\[\xymatrix{
w_*w^*f^*A_\lambda \ar[r] &
u_*u^*f^*A_\lambda \oplus
v_*v^*f^*A_\lambda \ar[r] &
j_*j^*f^*A_\lambda \ar[r] &
}\]
Taking the product over $\lambda$ and tensoring with $f^!E$ gives
a triangle, as does tensoring with $f^!E$ and then forming the product
over $\lambda$. There is a map between the triangles: if two
of the morphisms are isomorphisms then so is the third. Starting with
the fact that we know the map to be an isomorphism if $U\subset X$
is a sufficiently
small open affine, we glue to discover that it is an isomorphism for
every open immersion $u:U\la X$. The case of the identity map $\id:X\la X$
gives the Lemma.
\eprf

\lem{L99.3}
As in Lemma~\ref{L2.2}
let $U\stackrel\alpha\la V\stackrel\beta\la W$ be
finite-type morphisms of
noetherian schemes
so that $\alpha$ is a closed immersion and $\beta\alpha$
proper. Given a set of objects $\{A_\lambda,\,\lambda\in\Lambda\}$ in
the category $\Dqc(W)$, then the functors $\alpha^*$ and
$\alpha^\times$ take the natural morphism
\[\xymatrix@C+20pt{
\ds\beta^!\left[\prod_{\lambda\in\Lambda}A_\lambda\right]\ar[r]^-{J} &
\ds\prod_{\lambda\in\Lambda}\beta^! A_\lambda
}\]
to isomorphisms.
If $\beta$ is of finite Tor-dimension the functors
$\alpha^*$ and
$\alpha^\times$ also take the natural morphism
\[\xymatrix@C+20pt{
\ds\beta^*\left[\prod_{\lambda\in\Lambda}A_\lambda\right]\ar[r]^-{K} &
\ds\prod_{\lambda\in\Lambda}\beta^* A_\lambda
}\]
to isomorphisms.
\elem

\prf
Consider the commutative square
\[\xymatrix@C+20pt{
\ds\beta^\times
\left[\prod_{\lambda\in\Lambda}A_\lambda\right]\ar[r]^-{I}\ar[d]_-{\psi(\beta)} &
\ds\prod_{\lambda\in\Lambda}\beta^\times A_\lambda
\ar[d]^-{\prod_\Lambda \psi(\beta)}\\
\ds\beta^!\left[\prod_{\lambda\in\Lambda}A_\lambda\right]\ar[r]^-{J} &
\ds\prod_{\lambda\in\Lambda}\beta^! A_\lambda
}\]
The map $I$ is an isomorphism, after all $\beta^\times$ is a right
adjoint and respects products. If we apply the functor $\alpha^\times$
and recall (i) that it respects products and (ii) that
$\alpha^\times\psi(\beta)$ is an isomorphism by Lemma~\ref{L2.2}, then
we deduce the commutative diagram where the indicated maps are
isomorphisms
\[\xymatrix@C+20pt{
\ds \alpha^\times\beta^\times
\left[\prod_{\lambda\in\Lambda}A_\lambda\right]\ar[r]^-{\sim}\ar[d]_-{\wr} &
\ds \alpha^\times\left[\prod_{\lambda\in\Lambda}\beta^\times A_\lambda\right]
\ar[d]^-{\wr}\\
\ds \alpha^\times\beta^!\left[\prod_{\lambda\in\Lambda}A_\lambda\right]\ar[r]^-{\alpha^\times J} &
\ds \alpha^\times\left[\prod_{\lambda\in\Lambda}\beta^! A_\lambda\right]
}\]
Hence $\alpha^\times J$ must be an isomorphism. Support theory, more
concretely \cite[Proposition~A.3(ii)]{Iyengar-Lipman-Neeman13}, tells us that
$\alpha^*J$ is also an isomorphism. 

Now assume that $\beta$ is of finite Tor-dimension.
We wish to show that the maps $\alpha^*K$ and $\alpha^\times K$ are
isomorphisms, with $K$ as in the Lemma, and
by \cite[Proposition~A.3(ii)]{Iyengar-Lipman-Neeman13} it
suffices to consider $\alpha^*K$.

From
\ref{ST1.107.5}(ii) we have a natural isomorphism
$\beta^*(-)\oo\beta^!\co_W^{}\la\beta^!(-)$, hence the map
$J$ of the Lemma is isomorphic to the composite
\[\xymatrix@C+30pt{
\ds \beta^*\left[\prod_{\lambda\in\Lambda}A_\lambda\right]\oo\beta^!\co_W^{}\ar[r]^-{ K\oo\beta^!\co_W^{}} &
\ds \left[\prod_{\lambda\in\Lambda}\beta^* A_\lambda\right]\oo\beta^!\co_W^{}
\ar[r]^-{H} &
\ds \prod_{\lambda\in\Lambda}\left[\beta^* A_\lambda\oo\beta^!\co_W^{}\right]
}\]
Lemma~\ref{L99.-1} tells us that the map $H$ is an isomorphism, but
by the first part of the Lemma, which we already proved,
we know that $\alpha^*$ takes the composite to an
isomorphism. Hence $\alpha^*$ must take the map $K\oo\beta^!\co_W^{}$ to
an isomorphism, or to put it
differently $(-)\oo\alpha^*\beta^!\co_W^{}$ takes
$\alpha^*K$ to an isomorphism. We wish to show that
$\alpha^*K$ is an isomorphism and support theory,
more precisely \cite[Proposition~A.3(ii)]{Iyengar-Lipman-Neeman13},
tells us that
is suffices to show that every point of $U$ lies in the support
of $\alpha^*\beta^!\co_W^{}$. Since the support of $\alpha^*E$
is $\alpha^{-1}\supp(E)$ it certainly suffices
to show that $\supp(\beta^!\co_W^{})=V$.
But $\beta$ is of finite Tor-dimension, and
Lemma~\ref{L1.13}(iii) tells us that the support of $\beta^!\co_W^{}$
is equal to the support of $\beta^*\co_W^{}=\co_V^{}$, which is all of $V$.
\eprf

\lem{L99.5}
Let $U\stackrel\alpha\la V\stackrel\beta\la W$ be finite-type morphisms of
noetherian schemes
so that $\alpha$ is a closed immersion, $\beta$ is
of finite Tor-dimension, and $\beta\alpha$
proper. Suppose $E,F\in\Dqc(W)$ are any objects.
Then $\alpha^*$ and $\alpha^\times$ take the natural
map $\beta^*\HHom(E,F)\la\HHom(\beta^*E,\beta^*F)$ to an
isomorphism.
\elem

\prf
Support theory, more
concretely \cite[Proposition~A.3(ii)]{Iyengar-Lipman-Neeman13}, tells us that
$\alpha^*$ will take the map to an isomorphism if and only if
$\alpha^\times$ does. It suffices to prove the assertion for $\alpha^\times$.

Fix an object $F\in\Dqc(W)$ and let $\cl$ be the full subcategory of
all objects $E\in\Dqc(W)$ such that $\alpha^\times$ takes the map 
$\p_E^{}:\beta^*\HHom(E,F)\la\HHom(\beta^*E,\beta^*F)$ to an
isomorphism. If $E$ is a perfect complex then the map $\p_E^{}$ is an
isomorphism as it stands, hence all perfect complexes belong to $\cl$.
Also $\cl$ is obviously triangulated. Because $\Dqc(W)$ is compactly
generated
it suffices to prove that $\cl$ is localizing; we need to show that the
coproduct
of any set of objects in $\cl$ lies in $\cl$.

Assume therefore that $\{E_\lambda,\,\lambda\in\Lambda\}$ is a set
of objects of $\cl$. We wish to study the map
\[\xymatrix@C+30pt{
 \ds\beta^*\HHom\left[\left(\coprod_{\lambda\in\Lambda}
   E_\lambda\right)\,,\,F\right]\ar[r]^-{\p_{\coprod E_\lambda}^{}} &
  \ds\HHom\left[\beta^*\left(\coprod_{\lambda\in\Lambda}E_\lambda \right)\,,\,\beta^*F\right] 
}\]
which, up to isomorphism, rewrites as the composite
\[\xymatrix@C+30pt{
\ds\beta^*\prod_{\lambda\in\Lambda} \HHom(E_\lambda,F)\ar[r]^-{K}&
\ds\prod_{\lambda\in\Lambda} \beta^*\HHom(E_\lambda,F)\ar[r]^-{\prod_\lambda\p_{E_\lambda}^{}}&
\ds\prod_{\lambda\in\Lambda} \HHom(\beta^*E_\lambda, \beta^*F)
}\]
The fact that $\alpha^\times K$ is an isomorphism is 
by 
Lemma~\ref{L99.3}. The fact that $\alpha^\times
\big[\prod_\lambda\p_{E_\lambda}^{}\big]$ 
is an isomorphism
is because $\alpha^\times$ respects products, and takes each
$\p_{E_\lambda}^{}$ to an isomorphism.
\eprf

\section{Some more Hochschild-style formulas}
\label{S95}

To the extent that the
results surveyed in \S\ref{S1} and \S\ref{S2} are new, they
arose out of trying to
understand the magical formulas first discovered by
Avramov and Iyengar~\cite{Avramov-Iyengar08}, and developed further
in several papers, starting with Avramov, Iyengar, Lipman and
Nayak~\cite{Avramov-Iyengar-Lipman-Nayak10}. In this section we prove
some of these formulas---what
is novel is that they are stated and proved in the unbounded derived
category.
Since the results in the section are new---at least new in this
generality---we
present complete proofs. The reader interested in the highlights can
skip ahead to Corollary~\ref{C95.13} and Example~\ref{E95.15}.

\rmd{R95.1}
Let $f:X\la Y$ be a morphism of noetherian schemes. An object $M\in\Dqc(X)$
is called \emph{$f$--perfect} if 
\be
\item
$M$ belongs to $\dcoh(X)$; that is all but finitely many of the
cohomology sheaves vanish, and the non-vanishing ones are coherent.
\item
There exists an affine scheme $W$ and
a faithfully flat morphism $\rho:W\la X$ so that $(f\rho)_*\rho^*M$
is of finite Tor-amplitude. 
\ee
\ermd

\lem{L95.3}
Let $Z\stackrel i\la X\stackrel f\la Y$ be composable morphisms, with
$i$ a closed immersion
and $fi$ proper. Suppose that $M$ is an $f$--perfect object
in $\Dqc(X)$ and $C\in\Dqc(X)$ is a perfect complex supported on the
image of $i$. Then $f_*(M\oo C)$ is a perfect complex in $\Dqc(Y)$.
\elem

\prf
We need to show that $f_*(M\oo C)$ has coherent cohomology and is of
finite
Tor-amplitude. The coherence of the cohomology is clear: possibly
after replacing $Z$ by an infinitesimal thickening we can assume that
the complex $C\in\dcoh(X)$, whose cohomology is supported is on $Z$, is of the form
$i_*\wt C$ for some $\wt C\in\dcoh(Z)$. But then $M\oo C=M\oo i_*\wt
C\cong i_*(i^*M\oo \wt C)$, and hence $f_*(M\oo C)\cong f_*i_*(i^*M\oo
\wt C)$. But the complex $i^*M\oo \wt C\in\Dqcmi(Z)$ has coherent
cohomology,
and as $fi:Z\la Y$ is assumed proper so does $f_*i_*(i^*M\oo
\wt C)\cong f_*(M\oo C)$.

Now for the Tor-amplitude. Choose a faithfully flat map $\rho:W\la X$
as in Reminder~\ref{R95.1}. We know that there exists an integer $n$
so that the Tor-amplitude of $(f\rho)_*\rho^*M$ lies in the interval
$[-n,\infty)$, and therefore 
\[
(f\rho)_* \Big[(f\rho)^*\Dqc(Y)^{\geq 0}_{}\oo \rho^*M\Big]\eq 
\Dqc(Y)^{\geq 0}_{}\oo(f\rho)_*\rho^*M
\]
 is contained in
$\Dqc(Y)^{\geq -n}_{}$. The map $f\rho$ is a morphism from an affine
 scheme to $Y$, therefore
 $f\rho$ is quasi-affine. From the proof [not the statement]
 of \cite[Corollary~2.8]{Hall-Rydh13} we
deduce that $(f\rho)^*\Dqc(Y)^{\geq 0}_{}\oo
\rho^*M=\rho^*\big[f^*\Dqc(Y)^{\geq 0}_{}\oo M\big]$
is contained in $\Dqc(W)^{\geq -n}_{}$. Since $\rho^*$ is faithfully
flat it follows that $f^*\Dqc(Y)^{\geq 0}_{}\oo M$ is contained in
$\Dqc(X)^{\geq -n}_{}$.
As $C$ is assumed to be a perfect complex on $X$ its Tor-amplitude is
contained in the interval $[-m,m]$ for some $m>0$, and hence 
$f^*\Dqc(Y)^{\geq 0}_{}\oo M\oo C$ is contained in $\Dqc(X)^{\geq
  -m-n}_{}$.
But then
\[
f_*\Big[f^*\Dqc(Y)^{\geq 0}_{}\oo M\oo C\big] \eq \Dqc(Y)^{\geq
  0}_{}\oo f_*(M\oo C)
\]
is contained in $f_*\Dqc(X)^{\geq
  -m-n}_{}\subset \Dqc(Y)^{\geq
  -m-n}_{}$, and we deduce that the Tor-amplitude of $f_*(M\oo C)$ is
contained in $[-m-n,\infty)$.
\eprf

\rmd{R95.5}
Given any closed symmetric monoidal category $\cm$, that is a category $\cm$
with a symmetric tensor product and an internal Hom satisfying the usual
adjunction, there is a canonical evaluation map $\ev_{A,B}^{}:A\oo\HHom(A,B)\la
B$.
As is customary we deduce the following two 
canonical maps
\[\xymatrix@C+30pt@R-20pt{
\HHom(A,B)\oo C\ar[r]^-{\alpha} & \HHom(A,B\oo C) \\
A\oo\HHom(B,C) \ar[r]^-{\beta} & \HHom\big[\HHom(A,B)\,,\,C\big]
}\]
where $\alpha=\alpha(A,B,C)$ corresponds to the map
\[\xymatrix@C+30pt@R-20pt{
A\oo \HHom(A,B)\oo C\ar[r]^-{\ev_{A,B}^{}} & B\oo C 
}\]
and $\beta=\beta(A,B,C)$ corresponds to the composite
\[\xymatrix@C+30pt@R-20pt{
\HHom(A,B) \oo A\oo\HHom(B,C) \ar[r]^-{\ev_{A,B}^{}} &
B\oo\HHom(B,C)\ar[r]^-{\ev_{B,C}^{}} &C
}\]
An easy formal fact is that $\alpha(A,B,C)$ and $\beta(A,B,C)$ are isomorphisms
as long as $A$ is strongly dualizable.
\ermd

\con{C95.7}
We can of course use the basic maps of Reminder~\ref{R95.5} to
construct variants. The situation that interests us is where
$f:X\la Y$ is a morphism of schemes and $f^*\dashv f_*\dashv f^\times$
are the usual adjoint functors between $\Dqc(X)$ and $\Dqc(Y)$.
If $B,C\in\Dqc(Y)$ are any objects, then
the composite
\[\xymatrix@C+50pt@R-20pt{
  f^*\HHom(B,C)\oo f^\times B
  \ar[r]^-{\chi\big(\HHom(B,C),B\big)}&
  f^\times\big[\HHom(B,C)\oo B\big]
  \ar[r]^-{f^\times\ev_{B,C}^{}} &
f^\times C
}\]
induces by adjunction a map we
will denote $\gamma:f^*\HHom(B,C)\la\HHom(f^\times B,f^\times C)$.
If $A\in\Dqc(X)$ and $B,C\in\Dqc(Y)$ are objects, we can consider
the composites
\[\xymatrix@C+10pt@R-20pt{
  \HHom(A,f^\times B)\oo f^*C\ar[r]^-{\alpha} &
  \HHom(A, f^\times B\oo f^*C) \ar[r]^-{\HHom(A,\chi)} &
  \HHom\big[A,f^\times(B\oo C)\big]\\
  A\oo f^*\HHom(B,C)\ar[r]^-{\id\oo\gamma} &
  A\oo \HHom(f^\times B,f^\times C)\ar[r]^-{\beta} &
  \HHom\Big[\HHom(A,f^\times B)\,,\,f^\times C\Big]
}\]
\econ

\lem{L95.9}
Let $Z\stackrel i\la X\stackrel f\la Y$ be composable morphism of schemes
with $i$ a closed immersion and $fi$ a quasi-affine map. 
Let $A\in\Dqc(X)$ be an object,
and let $K\in\Dqc(X)$ be a perfect complex supported
on the closed subset $Z$, and so that $f_*(A\oo K)$ and $f_*(A\oo K^\vee)$
are perfect in $\Dqc(Y)$, where $K^\vee=\HHom(K,\co)$
is the dual of $K$. Then, with $B,C\in\Dqc(Y)$
arbitrary, the functor $K\oo(-)$ takes
the two composites at the end of Construction~\ref{C95.7}
to isomorphisms.
\elem

\prf
Because $K$ is perfect and supported on $Z$
it is isomorphic to a bounded complex
of coherent sheaves supported on $Z$. Up to replacing $Z$
by an infinitesimal thickening we may assume there exists an
object $\wt K\in\dcoh(Z)$ with $K\cong i_*\wt K$. But then we
have isomorphisms of functors
$K\oo(-)\cong i_*\wt K\oo(-)\cong i_*\big[\wt K\oo i^*(-)\big]$.
It certainly suffices to prove that $\wt K\oo i^*(-)$ takes
the maps at the end of Construction~\ref{C95.7} to isomorphisms. The
morphism $fi$ is quasi-affine and \cite[Corollary~2.8]{Hall-Rydh13}
tells us that $(fi)_*=f_*i_*$ is conservative, hence we are
reduced to showing that
$f_*i_*\big[\wt K\oo i^*(-)\big]\cong f_*\big[K\oo(-)\big]$
takes these two morphisms of Construction~\ref{C95.7}
to isomorphisms.

But now we are in business: tensoring these two maps with the strongly
dualizable $K$ has the effect of replacing $A$ with $A\oo K^\vee$ in the
first map and with $A\oo K$ in the second. Hence we are reduced to
showing that, under the assumption that $f_*A$ is a perfect complex
in $\Dqc(Y)$, the functor $f_*$ takes both of the  maps
of Construction~\ref{C95.7} to isomorphisms.
An easy exercise
with the standard isomorphisms $f_*(R\oo f^*S)\cong f_*R\oo S$
and $f_*\HHom(R,f^\times S)\cong \HHom(f_*R,S)$ allows
us to show that the functor $f_*$ takes the two maps
of Construction~\ref{C95.7}
to
\[\xymatrix@C+30pt@R-20pt{
\HHom(f_*A,B)\oo C\ar[r]^-{\alpha} & \HHom(f_*A,B\oo C) \\
f_*A\oo\HHom(B,C) \ar[r]^-{\beta} & \HHom\big[\HHom(f_*A,B)\,,\,C\big]
}\]
of Reminder~\ref{R95.5}, which are isomorphisms because $f_*A$ is perfect.
\eprf

Now we will apply our lemmas to the situation of Construction~\ref{C2.2002}.
We remind the reader:
 $f:X\la Y$ is a finite-type, flat morphism of noetherian schemes, 
and we formed the diagram
\[\xymatrix@C+10pt@R+10pt{
X \ar[dr]^-{\Delta} & & \\
 & X\times_Y^{}X \ar[r]^-{\pi_1^{}}\ar[d]_-{\pi_2^{}} & X\ar[d]^f \\
 & X \ar[r]^-{f} & Y\ar@{}[ul]|{(\diamondsuit)}
}\]
where the square is cartesian, $\pi_1^{}$ and $\pi_2^{}$ are the first
and second projections, and $\Delta:X\la X\times_Y^{}X$ is the
diagonal inclusion. We assert:

\pro{P95.11}
Let $A,B,C\in\Dqc(X)$ be objects and assume that $A$ is $f$--perfect.
Applying Construction~\ref{C95.7} to the morphism $\pi_2^{}:X\times_Y^{}X\la X$
and 
the objects $\pi_1^*A\in\Dqc(X\times_Y^{}X)$ and $B,C\in\Dqc(X)$,
we obtain morphisms
\[\xymatrix@C+50pt@R-20pt{
\HHom(\pi_1^*A,\pi_2^\times B)\oo \pi_2^*C\ar[r]^-{\HHom(\pi_1^*A,\chi)\circ\alpha} & \HHom\Big[\pi_1^*A,\pi_2^\times(B\oo C)\Big] \\
\pi_1^*A\oo\pi_2^*\HHom(B,C) \ar[r]^-{\beta\circ[\id\oo\gamma]} & \HHom\big[\HHom(\pi_1^*A,\pi_2^\times B)\,,\,\pi_2^\times C\big]
}\]
We assert that, if $L\in\Dqc(X\times_Y^{}X)$
is any object supported on the diagonal, then
the functors $L\oo(-)$ and $\HHom(L,-)$
take both maps to isomorphisms, as do the functors
$\Delta^*$ and $\Delta^\times$.
\epro

\prf
Let us first observe that the statement about $\Delta^*$ and $\Delta^\times$
follows from the assertion about $L\oo(-)$ and $\HHom(L,-)$. Since
$\Delta$ is a closed immersion the functor $\Delta_*$ is conservative,
and to show that $\Delta^*$ and $\Delta^\times$ take the
two maps to isomorphisms is equivalent to showing that
the composites $\Delta_*\Delta^*$ and $\Delta_*\Delta^\times$ take
them to isomorphisms. But it is standard that
$\Delta_*\Delta^*(-)\cong\Delta_*\co_X^{}\oo(-)$
 and
 $\Delta_*\Delta^\times(-)\cong\HHom(\Delta_*\co_X^{},-)$.
 Since $\Delta_*\co_X^{}$ is supported on the diagonal this reduces us
 to the statements about $L\oo(-)$ and $\HHom(L,-)$.

 Now let $\cl\subset\Dqc(X\times_Y^{}X)$ be the full subcategory of all
 objects $L$ so that $L\oo(-)$ and $\HHom(L,-)$ take both maps
 of the Proposition to isomorphisms. Clearly $\cl$ is a localizing
 subcategory. We wish to
 show that $\cl$ contains the category $\D_{\mathbf{qc},\Delta}^{}(X\times_Y^{}X)$,
that is the
full subcategory of $\Dqc(X\times_Y^{}X)$ of objects supported on the
diagonal. But the subcategory
$\D_{\mathbf{qc},\Delta}^{}(X\times_Y^{}X)$ is generated by
the objects inside it which are compact in the larger
$\Dqc(X\times_Y^{}X)$; this theorem was first proved in Thomason and
Trobaugh~\cite{ThomTro},
and for a more general, modern proof which works
for sufficienty nice algebraic stacks the reader can
see Hall and Rydh~\cite[Theorems~A, B and 4.10(2)]{Hall-Rydh13}.
This means
that any localizing subcategory, containing the compact
objects $K$ supported on the diagonal, will contain all
of $\D_{\mathbf{qc},\Delta}^{}(X\times_Y^{}X)$. It therefore suffices to
show that every compact $K$, supported on the
diagonal, belongs to $\cl$. Hence we let $K$ be a compact object
supported on the diagonal, and wish to show that $K\oo(-)$ and
$\HHom(K,-)\cong K^\vee\oo(-)$ take both maps in the
Proposition to isomorphisms.

Now the object $A\in\Dqc(X)$ is assumed $f$--perfect, and flat base-change
tells us that $\pi_1^*A$ is $\pi_2^{}$--perfect.
Consider the composable morphisms
$X\stackrel\Delta\la X\times_Y^{}X\stackrel{\pi_2^{}}\la X$; because
the composite $\id=\pi_2^{}\Delta$ is proper we may apply 
Lemma~\ref{L95.3}, and because it is quasi-affine Lemma~\ref{L95.9} also
applies. More precisely: with this pair of composable morphisms
apply Lemma~\ref{L95.3} to the $\pi_2^{}$--perfect object
$\pi_1^*A\in\Dqc(X\times_Y^{}X)$ and to the perfect complexes
$K,K^\vee\in\Dqc(X\times_Y^{}X)$
supported on the image of $\Delta$, and we learn that
$\pi_{2*}^{} (\pi_1^*A\oo K)$ and  $\pi_{2*}^{}(\pi_1^*A\oo K^\vee)$ are perfect in $\Dqc(X)$. But
then
Lemma~\ref{L95.9} allows us to conclude that $K\oo(-)$ and
$\HHom(K,-)\cong K^\vee\oo(-)$
take both morphisms of the Proposition to isomorphisms.
\eprf

\cor{C95.13}
Let $f:X\la Y$ be a finite-type, flat map of noetherian schemes, and
let the notation be as in Construction~\ref{C2.2002}.
For objects $A,C\in\Dqc(X)$ and $B\in\Dqc(Y)$, where $A$ is $f$--perfect, we have
isomorphisms
\begin{eqnarray*}
\HHom(A,f^!B)\oo C &\cong & \Delta^*\HHom\Big[\pi_1^*A,\pi_2^\times(f^*B\oo C)\Big] \\
\Delta^\times\Big[\pi_1^*A\oo\pi_2^*\HHom(f^*B,C)\Big]  &\cong& \HHom\big[\HHom(A,f^! B)\,,\,C\big]
\big]
\end{eqnarray*}
\ecor

\prf
The classical isomorphism
$\Delta^\times\HHom(E,F)\cong\HHom(\Delta^*E,\Delta^\times F)$,
coupled with the fact that $\Delta^\times\psi(\pi_2^{}):\Delta^\times\pi_2^\times\la\Delta^\times\pi_2^!$ is an
isomorphism by Lemma~\ref{L2.2}, tell us that for any
$E\in\Dqc(X\times_Y^{}X)$ and any $G\in\Dqc(X)$ the map
$\Delta^\times\HHom\big[E,\psi(\pi_2^{})\big]:\Delta^\times\HHom(E,\pi_2^\times
G)\la \Delta^\times\HHom(E,\pi_2^!G)$
is an isomorphism. By
\cite[Proposition~A.3(ii)]{Iyengar-Lipman-Neeman13}
we deduce that 
$\Delta^*\HHom\big[E,\psi(\pi_2^{})\big]:\Delta^*\HHom(E,\pi_2^\times
G)\la \Delta^*\HHom(E,\pi_2^!G)$
is also an isomorphism.

Now we turn to the proof of the Corollary. With the notation
of the Corollary, as a first step we prove
\begin{itemize}
\item
There is a natural isomorphism 
$\Delta^*\Big[\HHom(\pi_1^*A,\pi_2^\times f^*B)\Big]
\cong\HHom(A,f^!B)$.
\end{itemize}
The isomorphism of $\bullet$ comes from the following string of isomorphisms
\begin{eqnarray*}
\Delta^*\Big[\HHom(\pi_1^*A,\pi_2^\times f^*B)\Big]&\cong &
\Delta^*\Big[\HHom(\pi_1^*A,\pi_2^! f^*B)\Big]\\
&\cong&\Delta^*\HHom(\pi_1^*A,\pi_1^* f^!B)\\
&\cong&\Delta^*\pi_1^* \HHom(A,f^!B) \\
&\cong& \HHom(A,f^!B)\\
\end{eqnarray*}
The first isomorphism is because the
functor
$\Delta^*\HHom(\pi_1^*A,-)$ takes the map $\psi(\pi_2^{}):\pi_2^\times
f^*B\la \pi_2^!
f^*B$ 
to an isomorphism. The second isomorphism is because
$\theta(\diamondsuit):\pi_1^*f^!\la\pi_2^! f^*$ is an isomorphism,
see~\ref{ST1.107.9}(ii). The third isomorphism in by Lemma~\ref{L99.5},
and the last isomorphism is because $\Delta^*\pi_1^*\cong\id^*$.

With the preliminaries out of the way, apply
Proposition~\ref{P95.11} to the objects $A,f^*B,C\in\Dqc(X)$, where
$A$ is given to be
$f$--perfect.
The Proposition tells us that the functors $\Delta^*$ and
$\Delta^\times$ take the maps below to isomorphisms 
 \[\xymatrix@C+50pt@R-20pt{
\HHom(\pi_1^*A,\pi_2^\times f^*B)\oo \pi_2^*C\ar[r]^-{(1)} & \HHom\Big[\pi_1^*A,\pi_2^\times(f^*B\oo C)\Big] \\
\pi_1^*A\oo\pi_2^*\HHom(f^*B,C) \ar[r]^-{(2)} & \HHom\big[\HHom(\pi_1^*A,\pi_2^\times f^*B)\,,\,\pi_2^\times C\big]
}\]
And the first isomorphism of the Corollary is by applying $\Delta^*$
to the map $(1)$ while the second isomorphism is by applying
$\Delta^\times$ to the map $(2)$.  Let us take these one step at a
time, we begin be applying $\Delta^*$ to $(1)$. We obtain isomorphisms
\begin{eqnarray*}
\Delta^*\HHom\Big[\pi_1^*A,\pi_2^\times(f^*B\oo C)\Big] &\cong &
\Delta^*\Big[\HHom(\pi_1^*A,\pi_2^\times f^*B)\oo \pi_2^*C\Big] \\
&\cong&\Delta^*\Big[\HHom(\pi_1^*A,\pi_2^\times f^*B)\Big]\oo[\Delta^*\pi_2^* C]\\
&\cong& \HHom(A,f^!B)\oo C 
\end{eqnarray*}
The first isomorphism is just $\Delta^*$ applied to $(1)$. The second
isomorphism is because $\Delta^*$ respects the tensor product.
The third isomorphism is the tensor product of the isomorphism in $\bullet$
with the isomorphism
$\Delta^*\pi_2^*C\cong \id^*C=C$.

A similar analysis works for $\Delta^\times$ applied to the map
$(2)$, which gives us the first isomorphism below
\begin{eqnarray*}
\Delta^\times\Big[\pi_1^*A\oo\pi_2^*\HHom(f^*B,C) \Big] &\cong &\Delta^\times\HHom\big[\HHom(\pi_1^*A,\pi_2^\times f^*B)\,,\,\pi_2^\times C\big]\\
&\cong&\HHom\big[\Delta^*\HHom(\pi_1^*A,\pi_2^\times f^*B)\,,\,\Delta^\times\pi_2^\times C\big]\\
&\cong&\HHom\big[\HHom(A,f^!B)\,,\,C\big]
\end{eqnarray*}
The second isomorphism comes from the formula
$\Delta^\times\HHom(E,F)\cong\HHom(\Delta^*E,\Delta^\times F)$.
The third isomorphism is the functor $\HHom(-,-)$ applied
to the isomorphism in $\bullet$
and the isomorphism
$\Delta^\times\pi_2^\times C\cong \id^\times C=C$.
\eprf

\exm{E95.15}
Let us work out what Corollary~\ref{C95.13} says in the affine case:
that is $f:R\la S$ will be a finite-type, flat homomorphism of noetherian
rings and, by abuse of notation, we will also write $f:\spec S\la
\spec R$ for the induced map of noetherian schemes. 
We have the usual equivalences $\D(R)\cong\Dqc\big(\spec R\big)$ 
and $\D(S)\cong\Dqc\big(\spec S\big)$, and $f^*:\D(R)\la\D(S)$,
$f_*:\D(S)\la\D(R)$, $f^\times:\D(R)\la\D(S)$ and $f^!:\D(R)\la\D(S)$
are the affine versions of the standard functors of Grothendieck 
duality. Put $\se=S\oo_R^{} S$.
In this affine case, to say that an object  
$A\in\D(S)$ is $f$--perfect means that $A$ must have bounded
cohomology which is finite as $S$--modules,  and 
$f_*A\in\D(R)$ has finite Tor-dimension. 
 Let $A\in\D(S)$ be an
$f$--perfect complex, and let $B\in\D(R)$ and $C\in\D(S)$ be arbitrary.
Then the formulas of Corollary~\ref{C95.13} come down to
\begin{eqnarray*}
\Hom_S^{}(A,f^!B)\oo_S^{}C &\cong&S\oo_{\se}^{}\Hom_R^{}(A,B\oo_R^{}C)\ ,\\
\Hom_{\se}^{}\big[S,A\oo_R^{}\Hom_R^{}(B,C)\big] &\cong&
          \Hom_S^{}\big[\Hom_S^{}(A,f^!B),C\big]\ .
\end{eqnarray*}
where the Homs and tensors are all derived.
The reader can find special cases of these formulas
in Avramov, Iyengar, Lipman and Nayak~\cite{Avramov-Iyengar-Lipman-Nayak10}.

The reader might note that 
we have already met special cases of the
first of these formulas. If we put $A=C=S$
 then
the formula specializes to
\[
f^!B \quad\cong\quad
\Hom_S^{}(S,f^!B)\oo_S^{}S\quad\cong\quad S\oo_{\se}^{}\Hom_R^{}(S,S\oo_R^{}B)
\]
of the Introduction, and if we further specialize to $B=R$ we recover
the formula $f^!R\cong S\oo_{\se}^{}\Hom_R^{}(S,S)$
of Remark~\ref{R2.9}.
\eexm

\section{A historical review}
\label{S3}

Grothendieck first mentioned that he knew how to prove a relative version of
the Serre duality theorem in his ICM talk in Edinburgh in 1958, 
see~\cite{Grothendieck58C}. The first published version was
Hartshorne~\cite{Hartshorne66}; roughly speaking the construction
of $f^!$ given in \cite{Hartshorne66} is by gluing local data,
not an easy thing to do in the derived category. Three
and a half decades later 
Conrad~\cite{Conrad00} expanded and filled in details missing 
in~\cite{Hartshorne66}.
The presentation of the subject given here is entirely different in
spirit---it is based on early observations
by Deligne~\cite{Deligne66} and Verdier~\cite{Verdier68}, 
filled in and expanded greatly
in Lipman~\cite{Lipman09}. This second construction is much more global
and functorial, the usual objection to it is that it's difficult
to compute anything. 

Now it is time to say what's different here from the classical literature.
Let us begin with the observation that, until the late 1980s, no one
really understood how to handle unbounded derived categories. For the
first two decades of the subject the functor
$f^*$, which involves a derived tensor product, was treated as
a functor $f^*:\Dqcmi(Y)\la\Dqcmi(X)$, while the functor $f_*$, which 
involves injective resolutions, was classically viewed as a functor
$f_*:\Dqcpl(X)\la\Dqcpl(Y)$. A careful reader will note that, being defined
on different categories, these functors are not honest adjoints---there
is no counit of adjunction $f^*f_*\la\id$, and a classical version
of the treatment 
of \S\ref{S1} would have had to be more delicate. Luckily for us we live 
in modern times and can give the clean presentation of
the projection formula and the base-change maps of \S\ref{S1}. 

The article that brought modernity 
to this discipline was Spaltenstein~\cite{Spaltenstein88},
it taught us how to take injective and flat resolutions of unbounded complexes.
Spaltenstein's 
article made it clear how to define the adjoint functors $f^*:\Dqc(Y)\la\Dqc(X)$ and 
$f_*:\Dqc(X)\la\Dqc(Y)$. 
The natural question to arise was how much of Grothendieck duality could be
developed in the unbounded derived category. The existence of a right adjoint
$f^\times:\Dqc(Y)\la\Dqc(X)$ for $f_*$ 
was discovered soon after, the author even
showed in~\cite{Neeman96} that it is possible to obtain this adjoint 
easily and very formally using Brown representability. At the time the author 
was promoting the point of view
that the right way to approach all these classical
results was to employ systematically the techniques of homotopy theory,
like Brown representability---at the time this was still a novel idea. So
Lipman challenged the author to try to use the techniques of homotopy 
theory to extend Verdier's base-change theorem~\cite{Verdier68} to
the unbounded derived category. Instead of a proof
the author found a counterexample, 
see~\cite[Example~6.5]{Neeman96}. There exists a cartesian square
of noetherian schemes
\[
\xymatrix@C+10pt@R+10pt{
W \ar[r]^{u}\ar[d]_{f} & X\ar[d]^{g}\\
Y \ar[r]^{v} & Z\ar@{}[ul]|{(\diamondsuit)}
}
\]
with $v$ flat (even an open immersion) and $g$ proper
(even a closed immersion), and such that the base-change map
$\Phi(\diamondsuit):u^*g^\times\la f^\times v^*$ is \emph{not} an
isomorphism. As an aside we note that the schemes in question are all
affine. 

This counterexample had the unfortunate effect of
stifling the theory, for the next 
twenty years it put people off trying to develop the functor $f^!$ in the 
unbounded derived category. For example see Lipman's 
book~\cite{Lipman09}---Lipman makes a real effort to give the results
in the greatest generality in which they were known at the time, and 
for the functor $f^!$ he works
almost entirely with bounded-below complexes.
Drinfeld and Gaitsgory~\cite{Drinfeld-Gaitsgory13}
generalized a version of the theory to DG schemes,
and if the structure sheaf has negative cohomology then the category
$\Dqcpl(X)$ does not make much sense. To finesse the issue they work in
the category of Ind-coherent sheaves instead of the derived category.

In early 2013 I happened to run into Lipman at MSRI and he told me about
exciting recent work, joint with Avramov, Iyengar and Nayak, which
found a strange connection between Grothendieck duality and Hochschild 
homology and cohomology. In this survey we have already met this connection
in Theorem~\ref{T2.5} and
Example~\ref{E95.15}, see also Remarks~\ref{R2.7} and \ref{R2.9}.
 Theorem~\ref{T2.5} taught us about
this bizarre new map from Hochschild homology
to the dualizing complex $f^!\co_Y^{}$, and when $f$ is smooth and of relative
dimension $d$ this map happens to give an
isomorphism of $f^!\co_Y^{}$ with a  shift
of the relative canonical bundle. And in \S\ref{S95} we saw that
the formulas
of \S\ref{S2} are only the tip of the iceberg, there and many
more weird and wonderful ones---we presented two of them,
together with proofs, 
in Example~\ref{E95.15}.
The formulas of Example~\ref{E95.15}
are not new, special cases may be found 
in~\cite{Avramov-Iyengar08,Avramov-Iyengar-Lipman-Nayak10}.
What was new in \S\ref{S95} is that we gave them as special cases of
results that hold in the unbounded derived category.
Back in 2013,  when Lipman told me about the work, no one knew how
to define $f^!$ on the unbounded derived category.

Let us observe more carefully the second formula of Example~\ref{E95.15},
and for simplicity let's put $B=R$. 
We remind the reader, the formula is
\[
\Hom_{\se}^{}(S,A\oo_R^{}C) \quad\cong\quad
          \Hom_S^{}\big[\Hom_S^{}(A,f^!R),C\big]\ .
\]
If we fix $A$
and consider the expression on the right as a functor in $C$ then
it is clearly representable---the right hand side has the form
$\Hom_S^{}(P,-)$, where $P$ happens to be the expression $\Hom_S^{}(A,f^!R)$.
The isomorphism means that, as a functor in $C$, the expression 
$\Hom_{\se}^{}(S,A\oo_R^{}C)$ is also representable, in particular it commutes
with products---which is
far from obvious. The challenge Lipman gave me was to try to use Brown 
representability to prove these formulas.

There is such a proof, and Iyengar and I are working on writing it up. But 
this survey is about another direction our research took: in trying to
understand better these mysterious formulas we developed the natural
transformation $\psi(f):f^\times\la f^!$---early hints of it may be found
in Lipman~\cite[Exercise~4.2.3(d)]{Lipman09}.
What was new were the naturality and functoriality properties
of $\psi$, see~\cite{Iyengar-Lipman-Neeman13} for some illustrations
of their value. Because at the time $f^!$
was defined only on the bounded-below derived category our results imposed
artificial boundedness restrictions, and it was a natural challenge to try
to remove them. Working in the category of Ind-coherent sheaves, as in
Drinfeld and Gaitsgory, is clearly wrong for this problem---the formulas
of~\cite{Avramov-Iyengar-Lipman-Nayak10} live in the derived category.
The article \cite{Neeman13} was written to address this problem, in it
Grothendieck duality is developed in the unbounded derived category,
and we gave a brief summary of some of the results in \S\ref{S1}. In 
\S\ref{S2} and \S\ref{S95}
we gave illustrations of how one can approach the 
unbounded versions of the formulas
of~\cite{Avramov-Iyengar08,Avramov-Iyengar-Lipman-Nayak10,Iyengar-Lipman-Neeman13} using
the techniques surveyed in this paper---we
proved the formula $f^!=\Delta^*\pi_2^\times f^*$
for unbounded complexes in Proposition~\ref{P2.202},
while Example~\ref{E95.15} showed us how to derive the
reduction formulas of Avramov and Iyengar.
These formulas occur in
\cite{Avramov-Iyengar08,Avramov-Iyengar-Lipman-Nayak10,Iyengar-Lipman-Neeman13}, but with unnatural
boundedness hypotheses.

The reader might be puzzled. We mentioned that, twenty years ago, I produced
a counterexample~\cite[Example~6.5]{Neeman96} to the unbounded
version of Verdier's base-change theorem. There exists a cartesian square
of schemes
\[
\xymatrix@C+10pt@R+10pt{
W \ar[r]^{u}\ar[d]_{f} & X\ar[d]^{g}\\
Y \ar[r]^{v} & Z\ar@{}[ul]|{(\diamondsuit)}
}
\]
with $v$ an open immersion and $g$ proper, and such that the base-change map
$\Phi(\diamondsuit):u^*g^\times\la f^\times v^*$ is not an
isomorphism. So what has changed in two decades?
What's new is Theorem~\ref{T1.7}(ii): it tells us that, as long as
we further assume that $f$ is \emph{of finite Tor-dimension,} the problem goes
away and $\Phi(\diamondsuit)$ is an isomorphism. When we compare two
compactifications of $X$ we end up with cartesian squares of the
form
\[
\xymatrix@C+10pt@R+10pt{
X \ar[r]^{u}\ar[d]_{\id} & \ov X\ar[d]^{g}\\
X \ar[r]^{v} & \ov X'\ar@{}[ul]|{(\diamondsuit)}
}
\]
and the identity $\id:X\la X$ is of finite Tor-dimension. Thus the cartesian
squares that come up in the proof that $f^!=u^*p^\times$ is independent of
the factorization of $f:X\la Y$ as $X\stackrel u\la\ov X\stackrel p\la Z$
all have base-change maps which are isomorphisms.

The place where the old counterexample rears its ugly head is when it comes
to composition. The counterexample gave a commutative square (actually,
even cartesian) and hence we have $gu=vf$. Now $u$ and $v$ are open immersions
while $f$ and $g$ are proper, hence $u^!=u^*$, $v^!=v^*$, $f^!=f^\times$ and 
$g^!=g^\times$. On the other hand $u^!g^!=u^*g^\times$ is not isomorphic to
$f^!v^!=f^\times v^*$. We have already mentioned that, in
the old counterexample, the base-change map
$u^*g^\times\la f^\times v^*$ is not an isomorphism, but even more is true, the
functors are not isomorphic. The 2--functor $(-)^!$ is genuinely only 
oplax---there are natural maps $\rho(f,v):(vf)^!\la f^!v^!$ and 
$\rho(u,g):(gu)^!\la u^!g^!$, but clearly they cannot both be isomorphisms.
As it happens, in this
particular example $\rho(u,g)$ is an isomorphism while $\rho(f,v)$ isn't.

The situation is not hopeless: \ref{ST1.107.4} gives useful criteria for 
$\rho(f,g)$ to be an isomorphism, and in \S\ref{S2}, \S\ref{S99}
and \S\ref{S95}
we illustrated
how this can be applied to obtain unbounded versions of the results
of~\cite{Avramov-Iyengar08,Avramov-Iyengar-Lipman-Nayak10,Iyengar-Lipman-Neeman13}.
The illustrations of \S\ref{S2} also
showed how, with all these new methods, the abstract
nonsense approach to the subject pioneered by Deligne, Verdier and
Lipman can produce explicit computational formulas simply and easily. The
technicalities are different: the
``residual complexes'' of Grothendieck are replaced by the more
standard tools of Hochschild homology.

In passing let me note that
Hochschild homology is a {\it K}--theoretic invariant,
and its appearance raises the question whether more sophisticated 
{\it K}--theoretic invariants might shed even more light on the 
subject of Grothendieck duality. This is a volume on {\it K}--theory
and its applications to algebraic geometry, and it seems the 
appropriate place to raise this question.

\section{Generalizations}
\label{S4}

In March 2016 I received from the journal four referees' reports on my
article~\cite{Neeman13}. Mostly the referees' comments were simple enough
to address, but there were two difficult issues. One referee wondered what 
happens if we relax the noetherian hypothesis, while another suggested that 
I develop the entire theory in the generality of stacks.
This led me to think more carefully about these points. The noetherian
hypothesis seems indispensable, at least for this approach to the theory---some
of the lemmas have non-noetherian versions, but there is a point at which 
the argument runs into a brick wall without the noetherian assumption. 
But I'm happy to report that, under
relatively mild hypotheses, everything generalizes to noetherian 
algebraic stacks.

In fact the stacky version is cleaner to state. Algebraic stacks
naturally form a 2-category, as do triangulated categories.
The clean way to think about the theory is to view $(-)^*$, $(-)^\times$
and $(-)^!$ as 2-functors from [suitably restricted] algebraic stacks
to triangulated categories, with some relations among them. These relations
can be phrased in terms of the existence of certain natural transformations
relating these functors, and the assertion that certain
pairs of composites of natural transformations agree. For example:
it turns out that $(-)^*$ has the the structure of a monoid, meaning
there is a pseudonatural transformation $(-)^*\times(-)^*\la(-)^*$,
and $(-)^\times$ and $(-)^!$ are oplax modules over it. The
map $\psi$ turns out to be an oplax natural transformation
$\psi:(-)^\times\la(-)^!$ which is a module homomorphism. Anyway: the
reader can find a thorough discussion in the introduction to
\cite{Neeman13}.

This led to an expository conundrum in writing the current survey---it was unclear what was
the right generality for the results. Avramov, Iyengar, Lipman
and Nayak work with noetherian schemes, but allow the
morphisms to be essentially of finite type (rather than the more
restrictive finite type), and sometimes of finite Tor-dimension (rather than
flat). But the methods they use don't work for noetherian stacks---at
least not yet---because one doesn't yet know that a morphism of
noetherian stacks which is essentially of finite type has a Nagata
compactification. Nayak~\cite{Nayak09} proved the existence of Nagata
compactifications for 
morphisms of \emph{noetherian schemes} essentially of finite type, but so far no
one has generalized this even to algebraic spaces. In other words: the
results in this paper generalize in more than one direction, and at
present I do not know a common generalization that covers everything
that can be proved by the methods.

The compromise I made was to present the arguments in the
intersection of the known cases, that is finite-type, flat maps of noetherian
schemes, and leave to the reader the various generalizations. But I
did make an effort to give proofs that are easily adaptable, so that
the extension to (for example) algebraic stacks is straightforward.

When I gave the talk at TIFR, which amounted to a brief summary of
this survey, Geisser, Kahn, Saito
and Weibel raised the question of what portion of
the ideas might be transferrable to the six-functor formalism. Because
shortly after giving the talk I received the referees' reports, with the
questions about the non-noetherian and stacky versions of Grothendieck
duality, I haven't yet had the opportunity to think about this
other question. The six-functor
situation is another place where one defines functors like $f^!$ using
good factorizations of $f$, and it is eminently sensible to ask if there
might be an analog of the fancy, unbounded version of the base-change 
theorem and of its consequences.

The question is natural enough and I would like to come back to it when 
I have more time. In the interim I record it for others to study.

\def\cprime{$'$}
\providecommand{\bysame}{\leavevmode\hbox to3em{\hrulefill}\thinspace}
\providecommand{\MR}{\relax\ifhmode\unskip\space\fi MR }
\providecommand{\MRhref}[2]{%
  \href{http://www.ams.org/mathscinet-getitem?mr=#1}{#2}
}
\providecommand{\href}[2]{#2}

\end{document}